\newcommand\SEKImasterusepackages{\usepackage[T1]{fontenc}}
\newcommand\SEKIusersusepackages
{%
\usepackage{graphicx}
\usepackage{named,makeidx}
\usepackage{amssymb}
\bibliographystyle{named}
\def\citep{\cite}
\def\citet#1{\citeauthor{#1} \shortcite{#1}}
\newcommand\startcite{{\raise.2ex\hbox{[}}}
\newcommand\stopcite {\raise.2ex\hbox{]}}
\newcommand\citehelper[1]{\startcite #1\stopcite}

\newcommand\makeacitetoftwo[2]
{\citeauthor{#1} \citehelper{\citeyear{#1}; \citeyear{#2}}}

\newcommand\makeacitetofthree[3]
{\citeauthor{#1} \citehelper{\citeyear{#1}; \citeyear{#2}; \citeyear{#3}}}

\input dina4.sty
\input headerhot
\mathcommand\ident[1]{\mathsf{#1}}
\newcommand\plussymbol  {\ident{+}}
\newcommand\minussymbol {\ident{-}}
\newcommand\dividesymbol{\ident{/}}
\newcommand\timessymbol {\ident{*}}

\newcommand\set     {\ident{set}}

\newcommand\naturalssymbol{\ident{naturals}}
\newcommand\gensymsymbol{\ident{gensym}}
\mathcommand\mbpsymbol{\ident{m\hspace{-0.055em}b\hspace{-0.045em}p}}

\newcommand\csymbol     {\ident c}
\newcommand\esymbol     {\ident e}
\newcommand\fsymbol     {\ident f}
\newcommand\gsymbol     {\ident g}
\newcommand\hsymbol     {\ident h}
\newcommand\ksymbol     {\ident k}
\newcommand\psymbol     {\ident p}
\newcommand\ssymbol     {\ident s}
\newcommand\Everysymbol {\ident{Every}}
\newcommand\Permsymbol {\ident{Perm}}
\newcommand\RExistssymbol{\ident{Rexists}}
\newcommand\invertsymbol{\ident{invert}}
\newcommand\invsymbol{\ident{inv}}
\newcommand\abssymbol   {\ident{abs}}
\newcommand\cnssymbol   {\ident{cons}}
\mathcommand\cnsindexsymbol[1]{\ident{cons}_{#1}}
\newcommand\carsymbol   {\ident{car}}

\newcommand\cdrsymbol   {\ident{cdr}}
\newcommand\lengthsymbol{\ident{length}}
\newcommand\sizesymbol{\ident{size}}
\newcommand\dlsymbol    {\ident{dl}}
\newcommand\dloncesymbol{\ident{delfirst}}
\newcommand\rcsymbol    {\ident{rc}}
\newcommand\brsymbol    {\ident{br}}
\newcommand\revtailsymbol{\ident{revtail}}
\newcommand\revsymbol{\ident{rev}}
\newcommand\appendsymbol {\ident{append}}
\newcommand\zeropredicatesymbol{\ident{zerop}}
\newcommand\eqsymbol        {\ident{eq}}
\newcommand\ifthensymbol    {\mbox{\ident{If{}Then}}}
\newcommand\ifthenelsesymbol{\mbox{\ident{If{}ThenElse}}}
\mathcommand\eqindexsymbol        [1]{\eqsymbol        _{#1}}
\mathcommand\ifthenindexsymbol    [1]{\ifthensymbol    _{#1}}
\mathcommand\ifthenelseindexsymbol[1]{\ifthenelsesymbol_{#1}}
\newcommand\orsymbol    {\ident{or}}
\newcommand\andsymbol   {\ident{and}}
\newcommand\leqsymbol   {\ident{leq}}
\newcommand\lessymbol   {\ident{less}}
\newcommand\lexlessymbol{\ident{lexless}}
\newcommand\lexlimlessymbol{\ident{lexlimless}}
\newcommand\lexsymbol   {\ident{lex}}
\newcommand\acksymbol   {\ident{ack}}
\newcommand\switchsymbol{\ident{switch}}
\newcommand\swatchsymbol{\ident{swatch}}
\newcommand\diveinssymbol{\ident{div1}}
\newcommand\divzweisymbol{\ident{div2}}
\newcommand\divrestsymbol{\ident{divrest}}
\newcommand\diveinstailsymbol{\ident{div1tail}}
\newcommand\divzweitailsymbol{\ident{div2tail}}
\newcommand\remsymbol{\ident{rem}}
\newcommand\divsymbol{\ident{div}}

\newcommand\turingmachinesymbol{\ident T}
\newcommand\terminatespsymbol  {\ident{terminatesp}}
\newcommand\statesymbol        {\ident{state}}
\newcommand\cmdsymbol          {\ident{cmd}}
\newcommand\nthsymbol          {\ident{nth}}
\newcommand\doublesymbol       {\ident{double}}

\newcommand\ppsymbol           {\ident{p}}
\newcommand\qpsymbol           {\ident{q}}
\newcommand\Epsymbol           {\ident{E}}
\newcommand\Ppsymbol           {\ident{P}}
\newcommand\Qpsymbol           {\ident{Q}}
\newcommand\Marriessymbol      {\ident{Marries}}
\newcommand\Lovessymbol        {\ident{Loves}}
\newcommand\StolenBysymbol     {\ident{StolenBy}}
\newcommand\Humansymbol        {\ident{Human}}
\newcommand\Evensymbol         {\ident{Even}}
\newcommand\Oddsymbol          {\ident{Odd}}
\newcommand\Primesymbol        {\ident{Prime}}
\newcommand\EveryPairsymbol   {\ident{EveryPair}}
\newcommand\Givesymbol         {\ident{Give}}
\newcommand\Fathersymbol       {\ident{Father}}
\newcommand\Elephantpsymbol    {\ident{Elephant}}
\newcommand\Flowerpsymbol    {\ident{Flower}}
\newcommand\Germanpsymbol      {\ident{German}}
\newcommand\Bicyclepsymbol     {\ident{Bicycle}}
\newcommand\Hugepsymbol        {\ident{Huge}}
\newcommand\Animalpsymbol      {\ident{Animal}}
\newcommand\Malepsymbol        {\ident{Male}}
\newcommand\Boypsymbol         {\ident{Boy}}
\newcommand\Girlpsymbol        {\ident{Girl}}
\newcommand\Femalepsymbol      {\ident{Female}}
\newcommand\Roundpsymbol       {\ident{Round}}
\newcommand\Quadrangularpsymbol{\ident{Quadrangular}}
\newcommand\Metpsymbol         {\ident{Met}}
\newcommand\Kissedpsymbol      {\ident{Kissed}}
\newcommand\Bishopsymbol       {\ident{Bishop}}
\newcommand\mindexsymbol[1]{\existsvari w{#1}}

\newcommand\nonnegpsymbol      {\ident{nonnegp}}
\newcommand\wellsymbol         {\ident{well}}
\newcommand\welltailsymbol     {\ident{welltail}}
\newcommand\varsymbol          {\ident{var}}
\newcommand\aritysymbol        {\ident{arity}}

\newcommand\whilesymbol        {\ident{while}}

\newcommand\nullsymbol         {\ident{null}}
\newcommand\hdsymbol           {\ident{hd}}
\newcommand\tlsymbol           {\ident{tl}}
\newcommand\insymbol           {\ident{in}}
\newcommand\applysymbol        {\ident{app}}
\newcommand\termsymbol         {\ident{term}}
\newcommand\russellsymbol      {\ident{russell}}
\newcommand\sqrtindordsymbol[1]{\ident{sqrtio#1}}
\mathcommand\tightim{\longrightarrow}
\mathcommand\im{\ \tightim\ }
\mathcommand\rs{\:\rulesugar\:\:}
\mathcommand\rulesugar{\longleftarrow}

\mathcommand\doublepp[1]      {\doublesymbol   \beginargs{#1}\allargs}
\mathcommand\aritypp[1]      {\aritysymbol   \beginargs{#1}\allargs}
\mathcommand\lengthpp[1]      {\lengthsymbol   \beginargs{#1}\allargs}
\mathcommand\sizepp[1]      {\sizesymbol   \beginargs{#1}\allargs}
\mathcommand\wellpp[1]      {\wellsymbol   \beginargs{#1}\allargs}
\mathcommand\welltailpp[1]      {\welltailsymbol   \beginargs{#1}\allargs}
\mathcommand\varpp[1]      {\varsymbol   \beginargs{#1}\allargs}
\mathcommand\rempp[2]    {\remsymbol\beginargs{#1}\separgs{#2}\allargs}
\mathcommand\divpp[2]    {\divsymbol\beginargs{#1}\separgs{#2}\allargs}
\mathcommand\divrestpp[2]    {\divrestsymbol\beginargs{#1}\separgs{#2}\allargs}
\mathcommand\diveinspp[2]    {\diveinssymbol\beginargs{#1}\separgs{#2}\allargs}
\mathcommand\divzweipp[3]    {\divzweisymbol\beginargs{#1}\separgs{#2}
\separgs{#3}\allargs}
\mathcommand\diveinstailpp[4]    {\diveinstailsymbol\beginargs{#1}\separgs{#2}
\separgs{#3}\separgs{#4}\allargs}
\mathcommand\divzweitailpp[6]    {\divzweitailsymbol\beginargs{#1}\separgs{#2}
\separgs{#3}\separgs{#4}\separgs{#5}\separgs{#6}\allargs}
\mathcommand\mbppp[2]         {\mbpsymbol   \beginargs{#1}\separgs{#2}\allargs}
\mathcommand\revpp[1]     
{\revsymbol\beginargs{#1}\allargs}
\mathcommand\revppi[2]     
{\revsymbol^{#1}\beginargs{#2}\allargs}
\mathcommand\revtailpp[2]     
{\revtailsymbol\beginargs{#1}\separgs{#2}\allargs}
\mathcommand\revtailppi[3]
{\revtailsymbol^{#1}\beginargs{#2}\separgs{#3}\allargs}
\mathcommand\Permpp[2]     
{\Permsymbol\beginargs{#1}\separgs{#2}\allargs}
\mathcommand\Permppi[3]
{\Permsymbol^{#1}\beginargs{#2}\separgs{#3}\allargs}
\mathcommand\appendpp[2]      
{\appendsymbol \beginargs{#1}\separgs{#2}\allargs}
\mathcommand\appendppi[3]      
{\appendsymbol^{#1}\beginargs{#2}\separgs{#3}\allargs}
\mathcommand\Everypp[2]      
{\Everysymbol \beginargs{#1}\separgs{#2}\allargs}
\mathcommand\RExistspp[1]      
{\RExistssymbol \beginargs{#1}\allargs}
\mathcommand\appendlongpp[2]      
{\appendsymbol\left(\begin{array}{@{}l@{}}{#1}\separgs\\{#2}\end{array}\right)}
\mathcommand\cnspp[2]         {\cnssymbol   \beginargs{#1}\separgs{#2}\allargs}
\mathcommand\cnsppi[3]       {\cnssymbol^{#1}\beginargs{#2}\separgs{#3}\allargs}
\mathcommand\cnsindexpp[3]
{\cnsindexsymbol{#1}\beginargs{#2}\separgs{#3}\allargs}
\mathcommand\dlpp[2]          {\dlsymbol    \beginargs{#1}\separgs{#2}\allargs}
\mathcommand\dloncepp[2]      {\dloncesymbol\beginargs{#1}\separgs{#2}\allargs}
\mathcommand\dlonceppi[3]{\dloncesymbol^{#1}\beginargs{#2}\separgs{#3}\allargs}
\mathcommand\rcpp[2]          {\rcsymbol    \beginargs{#1}\separgs{#2}\allargs}
\mathcommand\brpp[2]          {\brsymbol    \beginargs{#1}\separgs{#2}\allargs}
\mathcommand\orpp[2]          {\orsymbol    \beginargs{#1}\separgs{#2}\allargs}
\mathcommand\andpp[2]         {\andsymbol   \beginargs{#1}\separgs{#2}\allargs}
\mathcommand\shortcnspp[2]    {\csymbol     \beginargs{#1}\separgs{#2}\allargs}
\mathcommand\tightshortcnspp[2]
{\csymbol\beginargs{#1}\tightsepargs{#2}\allargs}
\mathcommand\spp[1]           {\ssymbol     \beginargs{#1}\allargs}
\mathcommand\sppiterated[2]   {\ssymbol^{#1}\beginargs{#2}\allargs}
\mathcommand\sqrtindordpp[3]
                       {\sqrtindordsymbol{#1}\beginargs{#2}\separgs{#3}\allargs}
\mathcommand\ppp[1]           {\psymbol     \beginargs{#1}\allargs}
\mathcommand\pppiterated[2]   {\psymbol^{#1}\beginargs{#2}\allargs}
\mathcommand\zeropp           {\ident 0}
\mathcommand\Julietpp         {\ident{Juliet}}
\mathcommand\Romeopp          {\ident{Romeo}}
\mathcommand\Ipp              {\ident I}
\mathcommand\onepp            {\ident1}
\mathcommand\twopp            {\ident2}
\mathcommand\threepp          {\ident3}
\mathcommand\invertpp[1]      {\invertsymbol\beginargs{#1}\allargs}
\mathcommand\invpp[1]         {\invsymbol\beginargs{#1}\allargs}
\mathcommand\abspp[1]         {\abssymbol\beginargs{#1}\allargs}
\mathcommand\naturalspp[1]    {\naturalssymbol\beginargs{#1}\allargs}
\mathcommand\gensympp[1]      {\gensymsymbol\beginargs{#1}\allargs}
\mathcommand\nilpp            {\ident{nil}}
\mathcommand\falsepp          {\ident{false}}
\mathcommand\truepp           {\ident{true}}
\mathcommand\FALSEpp          {\ident{FALSE}}
\mathcommand\TRUEpp           {\ident{TRUE}}
\mathcommand\UNDEFpp          {\ident{UNDEF}}
\mathcommand\weirdppp         {\ident{weirdp}}
\mathcommand\ambigppp         {\ident{ambigp}}
\mathcommand\zeropredicatepp[1]{\zeropredicatesymbol\beginargs{#1}\allargs}
\mathcommand\cppeins       [1]{\csymbol     \beginargs{#1}\allargs}
\mathcommand\cppzwei       [2]{\csymbol\beginargs{#1}\separgs{#2}\allargs}
\mathcommand\eppeins       [1]{\esymbol     \beginargs{#1}\allargs}
\mathcommand\fppeins       [1]{\fsymbol     \beginargs{#1}\allargs}
\mathcommand\fppeinsindex  [2]{\fsymbol_{#1}\beginargs{#2}\allargs}
\mathcommand\fppeinsiterated[2]{\fsymbol^{#1}\beginargs{#2}\allargs}
\mathcommand\gppeins       [1]{\gsymbol     \beginargs{#1}\allargs}
\mathcommand\gppzwei       [2]{\gsymbol     \beginargs{#1}\separgs{#2}\allargs}
\mathcommand\hppeins       [1]{\hsymbol     \beginargs{#1}\allargs}
\mathcommand\kppeins       [1]{\ksymbol     \beginargs{#1}\allargs}
\mathcommand\appzero          {\ident a}
\mathcommand\bppzero          {\ident b}
\mathcommand\cppzero          {\ident c}
\mathcommand\dppzero          {\ident d}
\mathcommand\eppzero          {\ident e}
\mathcommand\eqindexpp[3]{\eqindexsymbol{#1}\beginargs{#2}\separgs{#3}\allargs}
\mathcommand\ifthenindexpp
[3]{\ifthenindexsymbol{#1}\beginargs{#2}\separgs{#3}\allargs}
\mathcommand\ifthenelseindexpp
[4]{\ifthenelseindexsymbol{#1}\beginargs{#2}\separgs{#3}\separgs{#4}\allargs}
\mathcommand\eqpp[2]{\eqsymbol\beginargs{#1}\separgs{#2}\allargs}
\mathcommand\leqpp[2]{\leqsymbol\beginargs{#1}\separgs{#2}\allargs}
\mathcommand\lespp[2]{\lessymbol\beginargs{#1}\separgs{#2}\allargs}
\mathcommand\lexlespp[2]{\lexlessymbol\beginargs{#1}\separgs{#2}\allargs}
\mathcommand\lexlimlespp[3]
               {\lexlimlessymbol\beginargs{#1}\separgs{#2}\separgs{#3}\allargs}
\mathcommand\lexpp[3]{\lexsymbol\beginargs{#1}\separgs{#2}\separgs{#3}\allargs}
\mathcommand\ackpp[2]{\acksymbol\beginargs{#1}\separgs{#2}\allargs}
\mathcommand\switchpp[1]{\switchsymbol\beginargs{#1}\allargs}
\mathcommand\swatchpp[1]{\swatchsymbol\beginargs{#1}\allargs}
\mathcommand\whilepp[2]{\whilesymbol\beginargs{#1}\separgs{#2}\allargs}
\mathcommand\nullpp[1]{\nullsymbol\beginargs{#1}\allargs}
\mathcommand\nullppiterated[2]{\nullsymbol^{#1}\beginargs{#2}\allargs}
\mathcommand\hdpp[1]{\hdsymbol\beginargs{#1}\allargs}
\mathcommand\hdppiterated[2]{\hdsymbol^{#1}\beginargs{#2}\allargs}
\mathcommand\carpp[1]{\carsymbol\beginargs{#1}\allargs}
\mathcommand\cdrpp[1]{\cdrsymbol\beginargs{#1}\allargs}
\mathcommand\tlpp[1]{\tlsymbol\beginargs{#1}\allargs}
\mathcommand\tlppiterated[2]{\tlsymbol^{#1}\beginargs{#2}\allargs}
\mathcommand\inpp[2]{\insymbol\beginargs{#1}\separgs{#2}\allargs}
\mathcommand\inppiterated[3]{\insymbol^{#1}\beginargs{#2}\separgs{#3}\allargs}
\mathcommand\applypp[2]{\applysymbol\beginargs{#1}\separgs{#2}\allargs}
\mathcommand\applyppiterated
[3]{\applysymbol^{#1}\beginargs{#2}\separgs{#3}\allargs}
\mathcommand\termpp[2]{\termsymbol\beginargs{#1}\separgs{#2}\allargs}
\mathcommand\setpp[1]{\set\beginargs{#1}\allargs}
\mathcommand\russellpp[1]{\russellsymbol\beginargs{#1}\allargs}

\mathcommand\Tpp[6]{\turingmachinesymbol\beginargs{#1}\separgs{#2}\separgs
{#3}\separgs{#4}\separgs{#5}\separgs{#6}\allargs}
\mathcommand\Tppseven[7]{\turingmachinesymbol\beginargs{#1}\separgs{#2}\separgs
{#3}\separgs{#4}\separgs{#5}\separgs{#6}\separgs{#7}\allargs}
\mathcommand\foreverppp[6]{\ident{foreverp}\beginargs{#1}\separgs{#2}\separgs
{#3}\separgs{#4}\separgs{#5}\separgs{#6}\allargs}
\mathcommand\terminatesppp[6]{\terminatespsymbol\beginargs{#1}\separgs
{#2}\separgs{#3}\separgs{#4}\separgs{#5}\separgs{#6}\allargs}
\mathcommand\terminatespppone[1]{\terminatespsymbol \beginargs{#1}\allargs}
\mathcommand\statepp
[3]{\statesymbol\beginargs{#1}\separgs{#2}\separgs{#3}\allargs}
\mathcommand\tightstatepp
[3]{\statesymbol\beginargs{#1}\tightsepargs{#2}\tightsepargs{#3}\allargs}
\mathcommand\cmdpp
[3]{\cmdsymbol  \beginargs{#1}\separgs{#2}\separgs{#3}\allargs}
\mathcommand\tightcmdpp
[3]{\cmdsymbol  \beginargs{#1}\tightsepargs{#2}\tightsepargs{#3}\allargs}
\mathcommand\stoppp           {\ident{stop}}
\mathcommand\leftpp           {\ident{left}}
\mathcommand\rightpp          {\ident{right}}
\mathcommand\nthpp         [2]{\nthsymbol  \beginargs{#1}\separgs{#2}\allargs}
\mathcommand\pppp          [1]{\ppsymbol\beginargs{#1}            \allargs}
\mathcommand\qppp          [2]{\qpsymbol\beginargs{#1}\separgs{#2}\allargs}
\mathcommand\Eppp          [1]{\Epsymbol\beginargs{#1}            \allargs}
\mathcommand\Epppzwei      [2]{\Epsymbol\beginargs{#1}\separgs{#2}\allargs}
\mathcommand\Pppp          [1]{\Ppsymbol\beginargs{#1}            \allargs}
\mathcommand\Ppppeinsindex [2]{\Ppsymbol_{#1}\beginargs{#2}\allargs}
\mathcommand\Ppppdrei      
[3]{\Ppsymbol\beginargs{#1}\separgs{#2}\separgs{#3}\allargs}
\mathcommand\Ppppvier
[4]{\Ppsymbol\beginargs{#1}\separgs{#2}\separgs{#3}\separgs{#4}\allargs}
\mathcommand\Qppp          [2]{\Qpsymbol\beginargs{#1}\separgs{#2}\allargs}
\mathcommand\Qpppeins      [1]{\Qpsymbol\beginargs{#1}\allargs}
\mathcommand\Qpppeinsindex [2]{\Qpsymbol_{#1}\beginargs{#2}\allargs}
\mathcommand\Qpppdrei      
[3]{\Qpsymbol\beginargs{#1}\separgs{#2}\separgs{#3}\allargs}
\mathcommand\Fatherpp      [2]{\Fathersymbol\beginargs{#1}\separgs{#2}\allargs}
\mathcommand\Marriespp     [2]{\Marriessymbol\beginargs{#1}\separgs{#2}\allargs}
\mathcommand\Lovespp       [2]{\Lovessymbol\beginargs{#1}\separgs{#2}\allargs}
\mathcommand\StolenBypp    [2]
{\StolenBysymbol\beginargs{#1}\separgs{#2}\allargs}
\mathcommand\Humanpp       [1]{\Humansymbol\beginargs{#1}\allargs}
\mathcommand\Evenpp        [1]{\Evensymbol\beginargs{#1}\allargs}
\mathcommand\Evenppi       [2]{\Evensymbol^{#1}\beginargs{#2}\allargs}
\mathcommand\Oddpp         [1]{\Oddsymbol\beginargs{#1}\allargs}
\mathcommand\Primepp       [1]{\Primesymbol\beginargs{#1}\allargs}
\mathcommand\EveryPairpp  [2]{\EveryPairsymbol\beginargs{#1}\separgs
{#2}\allargs}
\mathcommand\mindexppeins  [2]{\mindexsymbol{#1}\beginargs{#2}\allargs}
\mathcommand\Givepp        [3]{\Givesymbol
\beginargs{#1}\separgs{#2}\separgs{#3}\allargs}
\mathcommand\mindexppzwei  [3]{\mindexsymbol
{#1}\beginargs{#2}\separgs{#3}\allargs}
\mathcommand\mindexppdrei  [4]{\mindexsymbol
{#1}\beginargs{#2}\separgs{#3}\separgs{#4}\allargs}

\mathcommand\nonnegppp     [1]{\nonnegpsymbol\beginargs{#1}\allargs}

\mathcommand\anonymouscsymbol{c}
\mathcommand\anonymouscindexsymbol[1]{\anonymouscsymbol_{#1}}
\mathcommand\anonymousfsymbol{f}
\mathcommand\anonymouscpp
[2]{\anonymouscsymbol\beginargs{#1}\separgs\ldots\separgs{#2}\allargs}
\mathcommand\anonymouscindexpp
[3]{\anonymouscindexsymbol{#1}\beginargs{#2}\separgs\ldots\separgs{#3}\allargs}
\mathcommand\anonymousfpp
[2]{\anonymousfsymbol\beginargs{#1}\separgs\ldots\separgs{#2}\allargs}
\mathcommand\coerceindexpp[3]{[#3]_{#1}^{#2}}

\mathcommand\Elephantppp    [1]{\Elephantpsymbol\beginargs{#1}\allargs}
\mathcommand\Flowerppp      [1]{\Flowerpsymbol  \beginargs{#1}\allargs}
\mathcommand\Bicycleppp     [1]{\Bicyclepsymbol \beginargs{#1}\allargs}
\mathcommand\Germanppp      [1]{\Germanpsymbol  \beginargs{#1}\allargs}
\mathcommand\Hugeppp        [1]{\Hugepsymbol    \beginargs{#1}\allargs}
\mathcommand\Animalppp      [1]{\Animalpsymbol  \beginargs{#1}\allargs}
\mathcommand\Maleppp        [1]{\Malepsymbol    \beginargs{#1}\allargs}
\mathcommand\Boyppp         [1]{\Boypsymbol     \beginargs{#1}\allargs}
\mathcommand\Girlppp        [1]{\Girlpsymbol    \beginargs{#1}\allargs}
\mathcommand\Femaleppp      [1]{\Femalepsymbol  \beginargs{#1}\allargs}
\mathcommand\Roundppp       [1]{\Roundpsymbol   \beginargs{#1}\allargs}
\mathcommand\Bishoppp       [1]{\Bishopsymbol   \beginargs{#1}\allargs}
\mathcommand\Quadrangularppp[1]{\Quadrangularpsymbol  \beginargs{#1}\allargs}
\mathcommand\Kissedppp[2]{\Kissedpsymbol\beginargs{#1}\separgs{#2}\allargs}
\mathcommand\Metppp[2]   {\Metpsymbol   \beginargs{#1}\separgs{#2}\allargs}

\newcommand\bound     {{\rm bound}}
\newcommand\free      {{\rm free}}

\mathcommand\Vtripleindex[3]{\V\!_{{#1},\,{#2},\,{#3}}}
\mathcommand\Vdoubleindex[2]{\V\!_{{#1},\,{#2}}}
\mathcommand\Vsingleindex[1]{\V\!_{{#1}}}

\mathcommand\Erel[1]{\Gammaoffont\!_{#1}}
\mathcommand\Urel[1]{\Deltaoffont_{#1}}

\mathcommand\theRprimefromstrongtoweak{
  \inparenthesesinlinetight{
     \domres\id{\Vwall\cup\Vsome\setminus\RAN\varsigma}
     \nottight{\nottight\uplus}
     \reverserelation\varsigma
  }
  \nottight{\circ}
  \ranres
    {\transclosureinline R}
    {\Vwall\cup\Vsome\setminus\RAN\varsigma}
  \nottight{\nottight{\nottight{\uplus}}}
  \Vsome\tighttimes\Vsall
}

\mathcommand\deltaminus{\delta^-}
\mathcommand\deltaplus{\delta^+}
\mathcommand\deltaplusplus{\delta^{+^+}}
\mathcommand\deltastar{\delta^*}
\mathcommand\deltastarstar{\delta^{*^*}}

\mathcommand\Vall     {\Vsingleindex\indexdelta         }
\mathcommand\Vwall    {\Vsingleindex\indexdeltaminu     }
\mathcommand\Vsall    {\Vsingleindex\indexdeltaplus     }
\mathcommand\Vgsome   {\Vsingleindex\indexgammaplus     }
\mathcommand\Vsome    {\Vsingleindex\indexgamma         }
\mathcommand\Vfree    {\Vsingleindex\indexfree          }
\mathcommand\Vbound   {\Vsingleindex\indexbound         }
\mathcommand\Vsomesall{\Vsingleindex\indexgammadeltaplus}

\mathapplycommand\VARall      {\VARsingleindex\indexdelta         }
\mathapplycommand\VARwall     {\VARsingleindex\indexdeltaminu     }
\mathapplycommand\VARsall     {\VARsingleindex\indexdeltaplus     }
\mathapplycommand\VARgsome    {\VARsingleindex\indexgammaplus     }
\mathapplycommand\VARsome     {\VARsingleindex\indexgamma         }
\mathapplycommand\VARfree     {\VARsingleindex\indexfree          }
\mathapplycommand\VARbound    {\VARsingleindex\indexbound         }
\mathapplycommand\VARsomesall {\VARsingleindex\indexgammadeltaplus}
\mathcommand\displayVARsall[1]{\VARsingleindex\indexdeltaplus
\!\!\!\:\left(\begin{array}{@{}c@{}}#1\end{array}\right)}

\mathcommand\rigidvari     [2]{#1_{#2}^\indexgammadeltaplus}
\mathcommand\existsvari    [2]{#1_{#2}^\indexgamma    }
\mathcommand\forallvari    [2]{#1_{#2}^\indexdelta    }
\mathcommand\freevari      [2]{#1_{#2}^\indexfree     }
\mathcommand\wforallvari   [2]{#1_{#2}^\indexdeltaminu}
\mathcommand\sforallvari   [2]{#1_{#2}^\indexdeltaplus}
\mathcommand\gexistsvari   [2]{#1_{#2}^\indexgammaplus}
\mathcommand\boundvari     [2]{#1_{#2}}
\mathcommand\vari          [2]{#1_{#2}}
\mathcommand\wforallvarilow[2]{#1_{#2}^
{\raisebox{-.82ex}{\math\indexdeltaminu}}}

\newcommand\indexhelper[1]{{\scriptscriptstyle#1\:\!\!}}
\newcommand\indexdeltaplus
{\indexhelper{\delta^{\raisebox{-.17ex}{\fvesf\hskip-0.14em +}}}}
\newcommand\indexdeltaminu
{\indexhelper{\delta^{\mbox{\fvesf\hskip-0.14em\rule[.2ex]{.7em}{.15ex}}}}}
\newcommand\indexgammaplus
{\indexhelper{\gamma^{\mbox{\fvesf\hskip-0.14em +}}}}
\newcommand\indexgammadeltaplus
{\indexhelper{\gamma\delta^{\raisebox{-.17ex}{\fvesf\hskip-0.14em +}}}}

\newcommand\indexdelta{\indexhelper\delta}
\newcommand\indexgamma{\indexhelper\gamma}
\newcommand\indexfree
{{\scriptscriptstyle\free}}
\newcommand\indexbound
{{\scriptscriptstyle\bound}}

\newcommand\Wellfsymb{\ident{Wellf}}
\mathapplycommand\Wellfpp{\Wellfsymb}

\mathcommand\beginargs{(}
\mathcommand\allargs  {)}
\mathcommand\separgs  {,\,}
\mathcommand\tightsepargs{,}

\mathcommand\minusppnoparentheses  [2]{{#1}\,\minussymbol\,{#2}}
\mathcommand\tightminusppnoparentheses  [2]{{#1}\minussymbol{#2}}
\mathcommand\divideppnoparentheses [2]{{#1}\,\dividesymbol\,{#2}}
\mathcommand\plusppnoparentheses   [2]{{#1}\,\plussymbol \,{#2}}
\mathcommand\plusppnoparenthesesi  [3]{{#2}\,\plussymbol^{#1}\,{#3}}
\mathcommand\tightplusppnoparentheses   [2]{{#1}\plussymbol{#2}}
\mathcommand\timesppnoparentheses  [2]{{#1}\,\timessymbol\,{#2}}
\mathcommand\undppnoparentheses    [2]{{#1}\und            {#2}}
\mathcommand\oderppnoparentheses   [2]{{#1}\oder           {#2}}
\mathcommand\impliesppnoparentheses[2]{{#1}\implies        {#2}}
\mathcommand\leqinfixppnoparentheses[2]{{#1}\,\tight\leq\,{#2}}
\mathcommand\geqinfixppnoparentheses[2]{{#1}\,\tight\geq\,{#2}}
\mathcommand\dividepp [2]{(\divideppnoparentheses {#1}{#2})}
\mathcommand\minuspp  [2]{(\minusppnoparentheses  {#1}{#2})}
\mathcommand\pluspp   [2]{(\plusppnoparentheses   {#1}{#2})}
\mathcommand\timespp  [2]{(\timesppnoparentheses  {#1}{#2})}
\mathcommand\undpp    [2]{(\undppnoparentheses    {#1}{#2})}
\mathcommand\oderpp   [2]{(\oderppnoparentheses   {#1}{#2})}
\mathcommand\impliespp[2]{(\impliesppnoparentheses{#1}{#2})}

\input headertheorem
}

\newcommand\SEKIREVIEWSTATUS{1}

\newcommand\SEKISUBEDITOR
{{\sc Erica Melis} 
\\FR Informatik, Universit\"at des Saarlandes, D--66123 Saarbr\"ucken, Germany
\\E-mail: {\tt melis@activemath.org}
\\WWW: \url{http://www.activemath.org/~melis}
}

\newcommand\SEKIREVIEWER
{\nedoname
\\The \wittgenstein\ Archive, 3 Anderson Court, Newnham Road, Cambridge, CB3 9EZ, England
\\E-mail: {\tt mn15@cam.ac.uk}
}

\newcommand\SEKIREFEREEPROCESS
{LPAR 2003, 
 10th International Conference on Logic for Programming,
 Artificial Intelligence, and Reasoning;
 full paper, 15 pp.
}

\newcommand\SEKIDOCUMENTCLASS{0}

\newcommand\SEKIYEAR{2015}

\newcommand\SEKINUMBEROFYEAR{01}

\newcommand\SEKITYPE{0}

\newcommand\SEKIPUBLISHEDTYPE{0}
\newcommand\SEKIPUBLISHEDAS
{Master's thesis, 
 FR Informatik, Universit\"at des Saarlandes, D--66123 Saarbr\"ucken, Germany}
\title{\herbrandsfundamentaltheorem
\\---~an encyclopedia article~---} 
\author{\wirthnamenoindex\\\mbox{}\\\Institute\\\emailcp}

\input seki-deckblatt-5
\usepackage{amssymb}

\usepackage[T1]{fontenc}

\mathchardef\Gammaoffont="7000
\mathchardef\Gamma="0100
\mathchardef\Deltaoffont="7001
\mathchardef\Delta="0101
\mathchardef\Thetaoffont="7002
\mathchardef\Theta="0102
\mathchardef\Lambdaoffont="7003
\mathchardef\Lambda="0103
\mathchardef\Xioffont="7004
\mathchardef\Xi="0104
\mathchardef\Pioffont="7005
\mathchardef\Pi="0105
\mathchardef\Sigmaoffont="7006
\mathchardef\Sigma="0106
\mathchardef\Upsilonoffont="7007
\mathchardef\Upsilon="0107
\mathchardef\Phioffont="7008
\mathchardef\Phi="0108
\mathchardef\Psioffont="7009
\mathchardef\Psi="0109
\mathchardef\Omegaoffont="700A
\mathchardef\Omega="010A
\mathchardef\itype="017B

\catcode`\@=11

\gdef\allowhyphens{\penalty\@M \hskip\z@skip}

\gdef\set@low@box#1{\setbox\tw@\hbox{,}\setbox\z@\hbox{#1}\dimen\z@\ht\z@
     \advance\dimen\z@ -\ht\tw@
     \setbox\z@\hbox{\lower\dimen\z@ \box\z@}\ht\z@\ht\tw@ \dp\z@\dp\tw@ }
\gdef\set@low@boxsingle#1{\setbox\tw@\hbox{\rm,}\setbox\z@\hbox{#1}\dimen\z@\ht\z@
     \advance\dimen\z@ -\ht\tw@
     \setbox\z@\hbox{\lower\dimen\z@ \box\z@}\ht\z@\ht\tw@ \dp\z@\dp\tw@ }

\gdef\@glqq{%
\ifhmode\edef\@SF{\spacefactor\the\spacefactor}%
\else\let\@SF\empty
\fi
\CheckFamily\font\fraknomath\ifSameFamily ``\relax
\else\CheckFamily\font\swab\ifSameFamily ``\relax
\else\leavevmode\set@low@box{''}\box\z@\kern-.04em\allowhyphens\@SF\relax
\fi\fi}
\gdef\glqq{\protect\@glqq\kern+.07em}
\gdef\@grqq{%
\ifhmode\edef\@SF{\spacefactor\the\spacefactor}%
\else\let\@SF\empty 
\fi 
\CheckFamily\font\fraknomath\ifSameFamily ''\relax
\else\CheckFamily\font\swab\ifSameFamily ''\relax
\else\kern+.07em``\kern.07em\@SF\relax
\fi\fi}
\gdef\grqq{\protect\@grqq}
\gdef\@glq{{\ifhmode \edef\@SF{\spacefactor\the\spacefactor}\else
     \let\@SF\empty \fi \leavevmode
     \set@low@boxsingle{'\/}\box\z@\kern-.04em\allowhyphens\@SF\relax}}
\gdef\glq{\protect\@glq\kern+.07em}
\gdef\@grq{\ifhmode \edef\@SF{\spacefactor\the\spacefactor}\else
     \let\@SF\empty \fi \kern-.0125em`\kern.07em\@SF\relax}
\gdef\grq{\protect\@grq}

\catcode`\@=12

\newcommand\closequotecommanospace{''\nolinebreak\hskip-0.23em,}
\newcommand\closequotecomma      {\closequotecommanospace\         \,}

\makeatletter
\if \@ptsize 0
   \newfont{\scriptscriptscriptgoth}{ygoth scaled 760}
   \newfont{\scriptscriptgoth}{ygoth scaled 833}
   \newfont{\scriptgoth}{ygoth scaled 912}
   \newfont{\gothnomath}{ygoth}
   \newfont{\Goth}{ygoth scaled \magstephalf}
   \newfont{\GOth}{ygoth scaled \magstep1}
   \newfont{\GOTh}{ygoth scaled \magstep2}
   \newfont{\GOTH}{ygoth scaled \magstep3}

   \newfont{\scriptscriptscriptswab}{yswab scaled 760}
   \newfont{\scriptscriptswab}{yswab scaled 833}
   \newfont{\scriptswab}{yswab scaled 912}
   \newfont{\swab}{yswab}
   \newfont{\Swab}{yswab scaled \magstephalf}
   \newfont{\SWab}{yswab scaled \magstep1}
   \newfont{\SWAb}{yswab scaled \magstep2}
   \newfont{\SWAB}{yswab scaled \magstep3}

   \newfont{\scriptscriptscriptfrak}{yfrak scaled 760}
   \newfont{\scriptscriptfrak}{yfrak scaled 833}
   \newfont{\scriptfrak}{yfrak scaled 912}
   \newfont{\fraknomath}{yfrak}
   \newfont{\Frak}{yfrak scaled \magstephalf}
   \newfont{\FRak}{yfrak scaled \magstep1}
   \newfont{\FRAk}{yfrak scaled \magstep2}
   \newfont{\FRAK}{yfrak scaled \magstep3}

   \newfont{\init}{yinit}
   \newfont{\Init}{yinit scaled \magstephalf}
   \newfont{\INit}{yinit scaled \magstep1}
   \newfont{\INIt}{yinit scaled \magstep2}
   \newfont{\INIT}{yinit scaled \magstep3}
\fi
\if \@ptsize 1
   \newfont{\scriptscriptscriptgoth}{ygoth scaled 833}
   \newfont{\scriptscriptgoth}{ygoth scaled 912}
   \newfont{\scriptgoth}{ygoth}
   \newfont{\gothnomath}{ygoth scaled \magstephalf}
   \newfont{\Goth}{ygoth scaled \magstep1}
   \newfont{\GOth}{ygoth scaled \magstep2}
   \newfont{\GOTh}{ygoth scaled \magstep3}
   \newfont{\GOTH}{ygoth scaled \magstep4}

   \newfont{\scriptscriptscriptswab}{yswab scaled 833}
   \newfont{\scriptscriptswab}{yswab scaled 912}
   \newfont{\scriptswab}{yswab}
   \newfont{\swab}{yswab scaled \magstephalf}
   \newfont{\Swab}{yswab scaled \magstep1}
   \newfont{\SWab}{yswab scaled \magstep2}
   \newfont{\SWAb}{yswab scaled \magstep3}
   \newfont{\SWAB}{yswab scaled \magstep4}

   \newfont{\scriptscriptscriptfrak}{yfrak scaled 833}
   \newfont{\scriptscriptfrak}{yfrak scaled 912}
   \newfont{\scriptfrak}{yfrak}
   \newfont{\fraknomath}{yfrak scaled \magstephalf}
   \newfont{\Frak}{yfrak scaled \magstep1}
   \newfont{\FRak}{yfrak scaled \magstep2}
   \newfont{\FRAk}{yfrak scaled \magstep3}
   \newfont{\FRAK}{yfrak scaled \magstep4}

   \newfont{\init}{yinit scaled \magstephalf}
   \newfont{\Init}{yinit scaled \magstep1}
   \newfont{\INit}{yinit scaled \magstep2}
   \newfont{\INIt}{yinit scaled \magstep3}
   \newfont{\INIT}{yinit scaled \magstep4}
\fi
\if \@ptsize 2
   \newfont{\scriptscriptscriptgoth}{ygoth scaled 912}
   \newfont{\scriptscriptgoth}{ygoth}
   \newfont{\scriptgoth}{ygoth scaled \magstephalf}
   \newfont{\gothnomath}{ygoth scaled \magstep1}
   \newfont{\Goth}{ygoth scaled \magstep2}
   \newfont{\GOth}{ygoth scaled \magstep3}
   \newfont{\GOTh}{ygoth scaled \magstep4}
   \newfont{\GOTH}{ygoth scaled \magstep5}

   \newfont{\scriptscriptscriptswab}{yswab scaled 912}
   \newfont{\scriptscriptswab}{yswab}
   \newfont{\scriptswab}{yswab scaled \magstephalf}
   \newfont{\swab}{yswab scaled \magstep1}
   \newfont{\Swab}{yswab scaled \magstep2}
   \newfont{\SWab}{yswab scaled \magstep3}
   \newfont{\SWAb}{yswab scaled \magstep4}
   \newfont{\SWAB}{yswab scaled \magstep5}

   \newfont{\scriptscriptscriptfrak}{yfrak scaled 833}
   \newfont{\scriptscriptfrak}{yfrak}
   \newfont{\scriptfrak}{yfrak scaled \magstephalf}
   \newfont{\fraknomath}{yfrak scaled \magstep1}
   \newfont{\Frak}{yfrak scaled \magstep2}
   \newfont{\FRak}{yfrak scaled \magstep3}
   \newfont{\FRAk}{yfrak scaled \magstep4}
   \newfont{\FRAK}{yfrak scaled \magstep5}

   \newfont{\init}{yinit scaled \magstep1}
   \newfont{\Init}{yinit scaled \magstep2}
   \newfont{\INit}{yinit scaled \magstep3}
   \newfont{\INIt}{yinit scaled \magstep4}
   \newfont{\INIT}{yinit scaled \magstep5}
\fi

\newcommand{\mscriptscriptfrak}      [1]{\mbox{\scriptscriptscriptfrak#1}}
\newcommand{\mscriptfrak}            [1]{\mbox{\scriptscriptscriptfrak#1}}

\newcommand{\mfrak}[1]{\mbox{\fraknomath#1}}

\makeatother

\newif\ifSameFamily
\def\CheckFamily#1#2{\GetFamilyName{#1}\ArgOne
        \GetFamilyName{#2}\ArgTwo
        \ifx\ArgOne\ArgTwo\SameFamilytrue\else\SameFamilyfalse\fi}
\def\GetFamilyName#1{\edef\Tempa{#1}\def\Tempb{#1}\ifx\Tempa\Tempb
        \edef\Tempa{\fontname#1}\fi
        \edef\Tempa{\Tempa\space}%
        \expandafter\iGetFamilyName\Tempa\\}
\def\iGetFamilyName#1 #2\\#3{\def#3{#1}}
\def\DefFontName#1#2{{\escapechar-1\expandafter\expandafter\expandafter
        \iDefFontName\expandafter{\csname#2\endcsname}%
        \xdef#1{\expandafter\string\Tempa}}}
\def\iDefFontName{\def\Tempa}

\newcommand\unprotectedae
{\CheckFamily\font\fraknomath\ifSameFamily *a\else
 \CheckFamily\font\swab\ifSameFamily\char'212\else\"a\fi\fi}
\newcommand\unprotectedoe
{\CheckFamily\font\fraknomath\ifSameFamily 
*o\else\CheckFamily\font\swab\ifSameFamily\char'232\else\"o\fi\fi}

\DefFontName\eccclarge{eccc1200}
\DefFontName\eccc{eccc1000}
\DefFontName\ecccsmall{eccc0900}
\DefFontName\ecccfootnotesize{eccc0800}

\newcommand\unprotectedes
{\CheckFamily\font\fraknomath\ifSameFamily\char'215\else
\CheckFamily\font\swab\ifSameFamily\char'215\else  
s\fi\fi}

\newcommand\unprotectedesi
{\CheckFamily\font\fraknomath\ifSameFamily\char'215\else
\CheckFamily\font\swab\ifSameFamily\char'215\else  
\mbox{s\hskip.04em}\fi\fi
\-\ignorespaces}

\newcommand\unprotectedmyparagraphsymbol
{\CheckFamily\font\fraknomath\ifSameFamily 
\char'244\else\CheckFamily\font\swab\ifSameFamily
\char'244\else\S\fi\fi}

\renewcommand\ae{\protect\unprotectedae}
\renewcommand\oe{\protect\unprotectedoe}

\newcommand\es  {\protect\unprotectedes}

\newcommand\esi {\protect\unprotectedesi}  

\newcommand\myparagraphsymbol{\protect\unprotectedmyparagraphsymbol}

\def\mathfrak#1{%
\mathchoice
{{\mfrak{#1}}}
{{\mfrak{#1}}}
{{\mscriptfrak{#1}}}
{{\mscriptscriptfrak{#1}}}
}

\newcommand\namefont{}

\usepackage{url}

\newcommand\majorfootroom{\raisebox{-1.9ex}{\rule{0ex}{.5ex}}}

\newcommand\footroom{\raisebox{-1.5ex}{\rule{0ex}{.5ex}}}

\newcommand\smallfootroom{\raisebox{-1.0ex}{\rule{0ex}{3ex}}}

\newcommand\majorheadroom{\rule{0ex}{3.2ex}}
\newcommand\headroom{\rule{0ex}{2.8ex}}

\newcommand\claus    {Clau\es}

\newcommand\irving   {Irving}

\newcommand\jacques  {{\namefont Jacque\es}}

\newcommand\jean     {Jean}

\newcommand\michaelwithouttrema{Micha\-el}

\newcommand\peter    {Peter}

\newcommand\anellisindex    {\index{Anellis, Irving H. (1946--2013)}}

\newcommand\anellis         {{\namefont Anelli\es}}
\newcommand\anellisname     {{\namefont\anellisindex\irving\ H. \anellis}}

\newcommand\anellislifetime {(1946--2013)}

\newcommand\boole           {{\namefont Boole}}

\newcommand\myboolean       {\boole an}

\newcommand\brouwersecondname{Brouwer} 
\newcommand\brouwerfirstname{L.\,\,E.\,\,J.}
\newcommand\brouwerindex    {\index{Brouwer, L. E. J. (1881--1966)}}

\newcommand\brouwername     
                {\brouwerindex{\namefont \brouwerfirstname\ \brouwersecondname}}

\newcommand\fermatbirthyear 
{160?}

\newcommand\Begriffsschrift {Begriff\esi schrift}

\newcommand\gentzen         {{\namefont Gentzen}}

\newcommand\gentzensHauptsatz{{\namefont\gentzen}'\es\ {\Hauptsatz}}

\newcommand\Hauptsatz{Haupt\-satz}

\newcommand\goedel          {{\namefont G\oe del}}

\newcommand\secondincompletenesstheorem
{second \incompletenesstheorem}
\newcommand\secondIncompletenessTheorem
{Second \IncompletenessTheorem}

\newcommand\firstincompletenesstheorem
{first \incompletenesstheorem}
\newcommand\firstIncompletenessTheorem
{First \IncompletenessTheorem}
\newcommand\incompletenesstheorem{incompleteness theorem}
\newcommand\IncompletenessTheorem{Incompleteness Theorem}

\newcommand\heijenoortindex {\index{Heijenoort, Jean van (1912--1986)}}

\newcommand\heijenoort      {{\namefont Heijen\-oort}}
\newcommand\vanheijenoort   {{\namefont van Heijen\-oort}}
\newcommand\heijenoortname  {\heijenoortindex{\namefont\jean\ 
                             \vanheijenoort}}

\newcommand\herbrandindex     {\index{Herbrand, Jacques (1908--1938)}}

\newcommand\herbrand        {{\namefont Herbrand}}
\newcommand\herbrandname    {\herbrandindex{\namefont \jacques\ \herbrand}}
\newcommand\herbranddeathyear{1931}
\newcommand\herbrandlifetime{(1908--\herbranddeathyear)}

\newcommand\herbrandPhDtitle
{Recherche\es\ sur la th\'eorie de la d\'e\-mon\-stra\-tion}

\newcommand\herbrandsfundamentaltheoremindex
                                       {\index{Herbrand's Fundamental Theorem}}

\newcommand\herbrandsfundamentaltheoremnoindex
                                             {\herbrand'\es\ \fundamentaltheorem}
\newcommand\herbrandsfundamentaltheorem
           {\herbrandsfundamentaltheoremindex\herbrandsfundamentaltheoremnoindex}
\newcommand\herbrandsfalselemmaindex{\index{Herbrand's ``False Lemma''}}
\newcommand\herbrandsfalselemma
                           {\herbrandsfalselemmaindex\herbrand's ``False Lemma''}

\newcommand\fundamentaltheorem{Fundamental Theorem}

\newcommand\hilbert         {\mbox{\namefont Hilbert}}

\newcommand\loewenheim      {{\namefont L\oe wen\-heim}}

\newcommand\loewenheimskolem{\loewenheim--\skolem}
\newcommand\loewenheimskolemtheorem{\index
                    {Loewenheim-Skolem Theorem@L{\oe}wenheim--Skolem Theorem}%
                                                     \loewenheimskolem\ Theorem}

\newcommand\nedoindex       {\index{Nedo, Michael (1940)}}

\newcommand\nedo            {{\namefont Nedo}}
\newcommand\nedoname        {\nedoindex{\namefont\michaelwithouttrema\ \nedo}}

\newcommand\skolem          {{\namefont Skolem}}

\newcommand\skolemization   {\skolem\-ization}
\newcommand\skolemizedform  {\skolem\-ized form}   
\newcommand\skolemizedForm  {\skolem\-ized Form}   

\newcommand\PM              {Principia Mathematica}

\newcommand\wirth           {{\namefont Wirth}}
\newcommand\wirthnamenoindex{{\namefont\claus-\peter\ \wirth}}

\newcommand\wittgenstein    {{\namefont Wittgenstein}}

\newcommand\FB   {FB}

\newcommand\FBautinfveryshort{\FB\ AI}

\newcommand\f    {\mbox{}{f.}}   

\newcommand\ifandonlyif{if \nolinebreak and \onlyif}

\newcommand\onlyif{only \nolinebreak if}

\newcommand\getittotheright[1]  
{\hfill\mbox{}\penalty 100\mbox{\ \,}\nolinebreak
\hfill\hfill\hfill\nolinebreak\mbox{#1}\ignorespaces}

\newcommand\role{r\^ole}

\newcommand\theo {Theorem}

\newcommand\WWW  {WWW}

\newcommand\cf   {cf.}

\newcommand\cfnlb{\cf\nolinebreak}

\newcommand\CS   {Computer \Sci}

\newcommand\eg   {e.g.}

\newcommand\sententialtautology{sentential tautology}

\newcommand\firstorder{first-order}

\newcommand\FirstOrder{First-Order}

\newcommand\ie   {i.e.}

\newcommand\udiff{\ if\ }

\newcommand\May  {May}
\newcommand\Mar  {March}

\def\note{Note}

\newcommand\p    {p.}
\newcommand\pp   {pp.}
\newcommand\PP[2]{\pp\,\ignorespaces#1--\ignorespaces#2}

\newcommand\PhD  {PhD}
\newcommand\PhDthesis{\PhD\ thesis}
\newcommand\PhDThesis{\PhD\ Thesis}

\newcommand\sect {\myparagraphsymbol} 

\newcommand\Sci  {Sci.}

\newcommand\singulary{singulary}

\newcommand\wrt  {w.r.t.}

\newcommand\litnoteref[1]{\note\,#1}
\newcommand\litremaref[1]{Remark\,#1}

\newcommand\littheoref[1]{\theo\,#1}
\newcommand\litsectref[1]{\sect\,#1} 

\newcommand\litchapref[1]{Chapter\,#1} 

\newcommand\Examplename{Ex\-am\-ple}

\newcommand\litlemmref[1]{Lem\-ma\,#1}

\newcommand\remaref[1]{\litremaref{\ref{#1}}}

\newcommand\lemmref[1]{\litlemmref{\ref{#1}}}

\newcommand\theoref[1]{\littheoref{\ref{#1}}}

\newcommand\nthpositioner[2]
{#1\raisebox{0.52ex}{\scriptsize\hspace{0.07em}#2}}

\newcommand\mthpositioner[2]
{\math{#1}\raisebox{0.52ex}{\scriptsize\hspace{0.07em}#2}}
\newcommand\modulointocountzero[2]
{\count1=#1
\count2=#2
\count0=\count1
\divide  \count0 by \count2
\multiply\count0 by-\count2
\advance \count0 by \count1}
\newcommand\absolutevalueintocountzero[1]
{\count0=#1
\ifnum\count0<0\multiply\count0 by -1\fi}
\newcommand\nthstring[1]
{\def\myargone{#1}\ifcat a\myargone th\else\nthstringnochar{#1}\fi}
\newcommand\nthstringnochar[1]
{\absolutevalueintocountzero{#1}%
\modulointocountzero{\count0}{100} 
\ifnum\count0>9\ifnum\count0<20 th\else\stupidnthstring\fi
                                  \else\stupidnthstring\fi}
\newcommand\stupidnthstring
{\modulointocountzero{\count0}{10}
\ifnum\count0=1 \hskip-0.2em st\else
\ifnum\count0=2 nd\else
\ifnum\count0=3 rd\else 
                th\fi\fi\fi}

\newcommand\writeascents
[1]{\count4=#1
\ifnum\count4<0 
-\multiply\count4 by -1\fi
\modulointocountzero{\count4}{10}
\divide\count4 by 10
\count3=\the\count0
\modulointocountzero{\count4}{10}
\divide\count4 by 10
\the\count4
.\the\count0
\the\count3
}

\newcommand\frenchnthstring[1]
{\def\myargone{#1}\ifcat a\myargone th\else\frenchnthstringnochar{#1}\fi}
\newcommand\frenchnthstringnochar[1]
{\absolutevalueintocountzero{#1}%
\modulointocountzero{\count0}{100} 
\ifnum\count0>9\ifnum\count0<20 th\else\frenchstupidnthstring\fi
                                  \else\frenchstupidnthstring\fi}
\newcommand\frenchstupidnthstring
{\modulointocountzero{\count0}{10}
\ifnum\count0=1 \hskip-0.2em re\else
\ifnum\count0=2 me\else
\ifnum\count0=3 rd\else 
                th\fi\fi\fi}

\newcommand\CLAM      {{\rm CL\kern-.36em\raise.39ex\hbox{\sc a}\kern-.15emM}}

\newcommand\TEXMACS   {{\sc T\kern-.1667em\lower.5ex\hbox{E}\kern-.125emX\kern-.1em\lower.5ex\hbox{\textsc{m\kern-.05ema\kern-.125emc\kern-.05ems}}}}

\newcommand\Wernigerode    {Werni\-gerode}

\newcommand\plzwernigerode{\mbox{38855}}

\def       \emailcp      {{\tt wirth@logic.at}}

\newcommand\Institutedept
{\FBautinfveryshort}
\newcommand\Instituteinst
{\mbox{Hochschule Harz}}
\newcommand\Instituteplac
{\plzwernigerode\,\Wernigerode}
\newcommand\Institutecoun{Germany}

\newcommand\Institute
{\Institutedept, \Instituteinst, \Instituteplac, \Institutecoun}

\newcommand\academicpress{Academic Press (\elsevier)}

\newcommand\elsevier{Elsevier}

\newcommand\newspaperreference[5]
{\def\nameofjournalpress{#2}#1, #4 #5, #3\if?\nameofjournalpress
\else, #2\fi}

\newcommand\dateinjournal[1]{}

\newcommand\journalreference[6]
{\def\nameofjournalpress{#2}#1\nolinebreak\hskip.2em%
\dateinjournal{(#3) }{\mbox{\bf #4}}, \PP{#5}{#6}\if?\nameofjournalpress
\else, #2\fi}

\newcommand\journalreferenceprintyear[6]
{\def\nameofjournalpress{#2}#1 
#4:#5--#6, #3%
\if?\nameofjournalpress
\else, #2\fi}

\newcommand\journalreferenceprintyearaspartofnumber[6]
{\def\nameofjournalpress{#2}#1 
{#4/#3}, \PP{#5}{#6}\if?\nameofjournalpress
\else, #2\fi}

\newcommand\jscname
{J. Symbolic Computation}

\newcommand\jscprintyear
{\journalreferenceprintyear{\jscname}\academicpress}

\newcommand\tcsname{Theoretical \CS}
\newcommand\tcsjournal
{\journalreference\tcsname\elsevier}
\newcommand\tcsjournalprintyear
{\journalreferenceprintyear\tcsname\elsevier}

\urldef\wirthkuhnurl\url
{http://wirth.bplaced.net/SEKI/welcome.html#SWP-2007-01}

\urldef\urlsrninetythreedashzerofive\url
{http://wirth.bplaced.net/SEKI/welcome.html#SR-93-05}

\urldef\urlsreightyeighttwelve\url
{http://wirth.bplaced.net/SEKI/welcome.html#SR-88-12}

\date{\small Searchable Online Edition,
\\[+.4ex]Submitted \May\,24, 2014,
\\Revision Accepted \Mar\,4, 2015}
\newcommand\termsofdepth[2]{\app{{\mathcal T}_{#1}}{#2}}
\mathcommand\termsofdepthnovars[1]{{\mathcal T}_{#1}}
\mathcommand\termsofdepthnovarsstar[1]{{\mathcal T}_{#1}^*}
\newcommand\propertyC{\mbox{Property\hskip.2em C}} 
\renewcommand\namefont{\sc}
\newcommand\qed
{\relax\ifmmode\hskip2em \Box\else\unskip\nobreak\hskip1em $\Box$\fi}
\newcommand\existentialoid
{\math\gamma}\newcommand\existentialoidquantifier
{\mbox{\math\gamma-quantifier}}\newcommand\Anexistentialoidquantifier
{A \nolinebreak\existentialoidquantifier}\newcommand\anexistentialoidquantifier
{a \existentialoidquantifier}\newcommand
\Ruleofexistentialoidquantificationfirstpart
{Generalized rule of}\newcommand\Ruleofexistentialoidquantificationsecondpart
{\mbox{\math\gamma-quantification}}\newcommand
\Ruleofexistentialoidquantification
{\Ruleofexistentialoidquantificationfirstpart\
 \Ruleofexistentialoidquantificationsecondpart}\newcommand
\ruleofexistentialoidquantification
{generalized rule of \mbox{\math\gamma-quantification}}\newcommand
\Ruleofsimplificationfirstpart
{Generalized rule of}\newcommand\Ruleofsimplificationsecondpart
{simplification}\newcommand\Ruleofsimplification
{\Ruleofsimplificationfirstpart\ \Ruleofsimplificationsecondpart}\newcommand
\Ruleofuniversaloidquantificationfirstpart
{Generalized rule of}\newcommand
\Ruleofuniversaloidquantificationsecondpart
{\math\delta-quantification}\newcommand\Ruleofuniversaloidquantification
{\Ruleofuniversaloidquantificationfirstpart\ 
 \Ruleofuniversaloidquantificationsecondpart}\newcommand\universaloid
{\math\delta}\newcommand\universaloidquantifier
{\mbox{\math\delta-quantifier}}\newcommand\auniversaloidquantifier
{a \universaloidquantifier}
\begin{document}
\makecover
\maketitle
\begin{abstract}%
\sloppy
\herbrandsfundamentaltheorem\ provides a constructive characterization 
of derivability in first-order predicate logic by means of sentential logic.

Sometimes it is simply
called ``\herbrand's Theorem\closequotecomma but the longer name is
preferable as there are other important
``\herbrand\ theorems'' 
and \herbrand\ himself called it\,{\em ``Th\'eor\`eme fondamental''}.

It was ranked by \citet{bernays-herbrand} as follows:
``In its proof-theoretic form,  
\herbrandindex\herbrand's Theorem can be seen
as the central theorem of predicate logic.
It expresses the relation of predicate logic to propositional logic
in a concise and felicitous form.'' 
And by 
\citet{heijenoort-logic-calculus-language}:
``Let me say simply, 
in conclusion,
that {\em\Begriffsschrift}~\cite{begriffsschrift},
\loewenheim's paper~\shortcite{loewenheim-1915},
and Chapter\,\,5 of \herbrand's thesis \shortcite{herbrand-PhD}
are the three cornerstones of modern logic.''

\herbrandsfundamentaltheorem\ occurs in \litchapref 5 of 
his \PhDthesis\ \shortcite{herbrand-PhD}
---~entitled {\em\herbrandPhDtitle}~--- 
submitted by \herbrandname\ \herbrandlifetime\ 
in\,1929 at the University of Paris.

\herbrandsfundamentaltheorem\ is, 
together with \goedel's incompleteness theorems
and \gentzensHauptsatz,
one of the most influential theorems
of modern logic.

Because of its complexity,
\herbrandsfundamentaltheorem\ is typically 
fouled up in textbooks beyond all recognition.
As we are convinced that 
there is still much more to learn for the future from this theorem than many 
logicians know,
we will focus on the true message and its practical impact.
This requires a certain amount of streamlining of \herbrand's work,
which will be compensated by some remarks on the actual historical facts.
\end{abstract}
\vfill\pagebreak

\tableofcontents

\vfill\pagebreak

\section{Informal Introduction}
\vfill
\subsection{Validity in Sentential and in First-Order Logic}

The language of classical (\ie\ two-valued)
{\em sentential logic} (also called ``propositional logic'')
is \nolinebreak formed by {\em\myboolean\ operator symbols}
\\\LINEnomath
{---~~~~say conjunction \nlbmaths\tightund, \
disjunction \nlbmaths\tightoder, \
negation \nlbmath\neg~~~~---}
\\on {\em sentential variables}
(\ie\ nullary predicate symbols). \hskip.3em
For simplicity, but without loss of generality, 
we \nolinebreak will consider exactly the these three operators symbols 
as part of our language of sentential logic in this
article. \hskip.2em
Other operators will be considered just as syntactical sugar; \hskip.2em
for instance, material implication
\bigmaths{A\tightimplies B}{} will be considered a 
meta-level notion defined as \bigmaths{\neg A\tightoder B}.
The interpretation of the \myboolean\ operator symbols is fixed,
whereas the sentential variables range 
over the \myboolean\ values \TRUEpp\ and \FALSEpp. \hskip.3em
A \nolinebreak sentential formula is {\em valid}\/ \hskip.1em
if it evaluates to \nlbmath\TRUEpp\ 
for \nolinebreak all interpretations 
(\ie\ mappings to \myboolean\ values) 
of the sentential variables.

In a first step,
let us now add non-nullary {\em predicate symbols}, 
which take terms as arguments. 
{\em Terms}\/ are formed from {\em function symbols} and {\em variables}\/
over a non-empty {\em domain of individuals}, 
which has to be chosen by any interpretation and 
is a assumed to be well-determined and fixed in advance,
although it may be infinite. \hskip.2em
Such a {\em quantifier-free first-order formula}\/ 
is {\em valid}\/ \hskip.1em
if it evaluates to \nlbmath\TRUEpp\
for all interpretations of predicate symbols as functions from individuals 
to \myboolean\ values, of function symbols as functions from individuals to 
individuals, and of variables as individuals.

Note that this extension is not a substantial one, however, 
because the notion of validity does not change when 
we interpret the quantifier-free first-order formulas as sentential formulas, 
simply by considering the predicates together with their argument terms just
as names for atomic sentential variables.

In a second step,
we can add quantifiers such as ``\math{\forall\,}'' (``for all \ldots'') and 
``\math{\exists\,}''
(``there is a \ldots'') to bind variables. 
This means that formulas are now formed not only by applying \myboolean\ operators
to formulas, 
but also the \singulary\ operators ``\math{\forall x.}'' and 
\math{\exists x.}'', 
{\em binding}\/ an arbitrary variable symbol \nlbmaths x. \hskip.4em
Evaluation is now defined for these additional formula formations 
in the obvious way: \hskip.3em
\bigmath{\exists x\stopq A} 
(or else: \math{\forall x\stopq A}) \hskip.1em
evaluates to \TRUEpp\ if 
the single formula argument \nlbmath A 
\mbox{(its {\em scope})} \hskip.2em
evaluates to \TRUEpp\
for some interpretation of \nlbmath x 
(or else: for \nolinebreak all interpretations of \nlbmath x); \hskip.4em
otherwise it evaluates to \nlbmaths\FALSEpp.

With this second step we arrive at {\em first-order predicate logic}\/
(with function symbols).
This logic is crucially different from sentential logic,
because the testing of all domains of individuals becomes now 
unavoidable for determining validity of a formula in general. \hskip.3em
Even though it \nolinebreak actually suffices 
to check only one domain for each cardinality
(different from \nlbmaths 0, but including infinite ones), \hskip.2em
this cannot be executed effectively in general.
As noted above, 
however,
the domains do not matter 
if no quantifiers occur in a first-order formula.
\begin{definition}[Sentential Validity]\\
A first-order formula is {\em sententially valid}
\udiff\ it is quantifier-free and valid in sentential logic,
provided that we consider 
the predicates together with their argument terms just
as names for atomic sentential variables.
\getittotheright\qed\end{definition}
Note that a formula does not change its meaning if we replace a bound
variable with a fresh one.
For instance, there is not difference in validity between
\\\LINEmaths{\forall x\stopq\inparentheses
{\app{\ident{Human}}x\implies\app{\ident{Mortal}}x}}{}
\\and
\LINEmaths{\forall y\stopq\inparentheses
{\app{\ident{Human}}y\implies\app{\ident{Mortal}}y}},\mbox{~~~~~}
\\both expressing that ``all humans are mortal'' 
---~in a structure 
where the \singulary\ predicates {\ident{Human}} and {\ident{Mortal}}
have the obviously intended interpretation.
Note, however, that none of these equivalent formulas is
valid, because we also have to consider the 
structure where {\ident{Human}} is always \TRUEpp\ and {\ident{Mortal}}
is \FALSEpp, in which case the formula evaluates to \nlbmaths\FALSEpp.

Just like \herbrand,
we consider equality of formulas only up to renaming of bound variables.
Thus, we consider the two displayed formulas to be identical.

\begin{sloppypar}
A variable may also occur {\em free}\/ in a formula,
\ie\ not in the scope of any quantifier binding it. \hskip.2em
We \nolinebreak will, 
however,
tacitly consider only formulas where each occurrence of 
each variable is either free or otherwise bound by a unique quantifier.
This excludes ugly formulas such \nolinebreak as 
\mathnlb{\app{\ident{Human}}x
\hskip.06em\tightund\hskip.11em\exists x.\app{\ident{Mortal}}x}{,\,}
\mathnlb{\exists x.\app{\ident{Human}}x
\hskip.06em\tightund\hskip.11em
\exists x.\app{\ident{Mortal}}x}, or
\mathnlb{\forall x.\inparenthesestight{\app{\ident{Human}}x
\hskip.06em\tightund\hskip.11em\exists x.\app{\ident{Mortal}}x}}.
The bound variables of such formulas can always be renamed to obtain
nicer formulas in our restricted sense,
such as \bigmathnlb
{\app{\ident{Human}}x\und\exists z.\app{\ident{Mortal}}z}, 
\bigmathnlb{\exists x.\app{\ident{Human}}x\und
\exists z.\app{\ident{Mortal}}z}, and
\bigmathnlb{\forall x.\inparenthesestight
{\app{\ident{Human}}x\und\exists z.\app{\ident{Mortal}}z}}.
Both human comprehension and formal treatment become less difficult 
by this common syntactical restriction.\end{sloppypar}

\subsection{Calculi: Soundness, Completeness, Decidability}

To get a more constructive access to first-order predicate logic,\
validity has to be replaced with derivability in a calculus.
Such a calculus is {\em sound}\/ if we can derive only valid formulas with it,
and {\em complete}\/ if every valid formula can be derived with it.
Luckily, there are sound and complete calculi for first-order logic.

Let us consider formal derivation in a sound and complete calculus 
for first-order logic.
Then there are effective enumeration procedures that, in the limit, 
would produce an infinite list
of all derivable consequences.
This means that derivability in first-order logic is {\em semi-decidable}\/:
If we want to find out whether a first-order formula is derivable, 
we can start such an enumeration procedure
and say ``yes'' if our formula comes along.

Non-derivability in first-order logic, however, is not semi-decidable:
There cannot be an 
enumeration procedure for those first-order formulas which are not derivable.
\hskip.2em
In \nolinebreak other words, 
derivability is not {\em co-semi-decidable}.

A problem is {\em decidable}\/ if it is both semi- and co-semi-decidable. 
Therefore, the problem of derivability in first-order logic 
is not decidable:
There cannot be any effective procedure that,
for an arbitrary first-order formula as input,
always returns that one of the answers 
``yes'' and ``no'' that is correct \wrt\ its derivability.
Note that the problem of derivability in first-order logic 
is historically called the 
{\em Entscheidung\esi problem (in engerer Bedeutung)}, 
\ie\ the {\em decision problem (for first-order logic)});
\cfnlb\ 
\cite[\p,72\f]{grundzuege},
\cite
[\litnoteref{8.6}, \p\,8]{grundlagen-german-english-edition-volume-one-one}.

Sentential logic, however,
is decidable.

Therefore, it makes sense to characterize derivability in first-order logic
by a semi-decision procedure based on validity in sentential logic.
\pagebreak

\begin{remark}[Historical Correctness]\\
The notion of decidability was developed mainly after \herbrand's death.
The {\em Ent\-schei\-dung\esi problem}\/ was an open
problem during \herbrand's lifetime,
because the co-semi-undecidability was established only later by 
\citet{Church36b} and 
\citet{turing-entscheidungsproblem}.
\getittotheright\qed\end{remark}

\subsection{First Major Aspect of \herbrandsfundamentaltheorem}
A major aspect of \herbrandsfundamentaltheorem\ is that it
provides a semi-decision procedure for first-order logic as follows:
For a given first-order formula \maths A,
this procedure produces
a list of quantifier-free first-order formulas
\par\noindent\LINEmaths{
F^{\,\app{\termsofdepthnovars 1}F}\comma
F^{\,\app{\termsofdepthnovars 2}F}\comma
F^{\,\app{\termsofdepthnovars 3}F}\comma
\ldots}{}
\\[+.3ex]\noindent such that \math A is derivable in first order-logic
\ifandonlyif\ one of the formulas 
\nlbmath{F^{\,\app{\termsofdepthnovars i}F}} 
is sententially valid. \hskip.3em
We say that 
\math A has {\em\propertyC\ of order \nlbmath i}
\udiff\ \math{F^{\,\app{\termsofdepthnovars i}F}} is sententially valid.

\subsection{Second Major Aspect of \herbrandsfundamentaltheorem}
Another major aspect of \herbrandsfundamentaltheorem\ is
that in {\em\herbrand's {\it modus ponens}-free calculus}\/ for first-order logic
there is a linear derivation of \math A from 
\nlbmath{F^{\,\app{\termsofdepthnovars i}F}},
provided that \math A has \propertyC\ of order \nlbmaths i. \
A derivation is {\em linear}\/ if
---~seen as a tree~--- 
it \nolinebreak has no branching
because all inference rules have exactly one premise. \hskip.2em
In addition, 
this \nolinebreak derivation also has the so-called ``sub''-formula 
property \wrt\ \nlbmath A. \hskip.3em
Moreover,
contrary to all calculi that were invented before,
and similar to the calculi of \cite{gentzen},
\hskip.2em
\herbrand's {\it modus ponens}\/-free calculus
gives humans a good chance to actually find this linear derivation 
based on an informal proof. \hskip.3em
Furthermore, 
\herbrand's {\it modus ponens}\/-free calculus shows a great similarity with
today's approaches to automated theorem proving,
greater even than that of the well-known calculi of \cite{gentzen}.


\subsection{Also a Completeness Theorem for \FirstOrder\ Logic}
``\propertyC\hskip.07em'' is a name introduced in \cite{herbrand-PhD}.
Without a name, 
this property occurs already in \cite{loewenheim-1915}, \hskip.2em
where it is shown that a first-order formula 
is valid
\ifandonlyif\ it has \propertyC\ of order \nlbmaths i, \hskip.2em
for some positive natural number \nlbmath i 
\ ---~~which became famous as the \loewenheimskolemtheorem.

In his \PhDthesis,
\herbrandindex\herbrand\ also showed the equivalence of 
his own first-order calculi with those of the \hilbert\ school
\cite{grundlagen-german-english-edition-volume-one-two}
and the \PM~\cite{PM}. \
Therefore,
as a consequence of the \loewenheimskolemtheorem,
the completeness of all these calculi is an immediate
corollary of \herbrandsfundamentaltheorem. \

\herbrandindex\herbrand,
however, 
did not trust the notion of first-order validity.
As the first follower of \hilbert's {\em finitistic 
standpoint}\/ in proof theory in France,
\herbrandindex\herbrand\ was so radically finitistic that 
---~in the area of logic~--- 
he did not accept model theory or set theory at \nolinebreak all. \
And so \goedel\ proved the completeness of first-order logic first 
when he
submitted his thesis~\cite{goedel-completeness} in 1929, 
in the same year as \herbrandindex\herbrand, and the
theorem is now called {\em \goedel's Completeness Theorem}\/ 
in all textbooks on logic.%
\pagebreak

\begin{figure}[h]
\begin{center}%
\includegraphics
[width=130mm,height=55mm,bb=0 0 459 238]%
{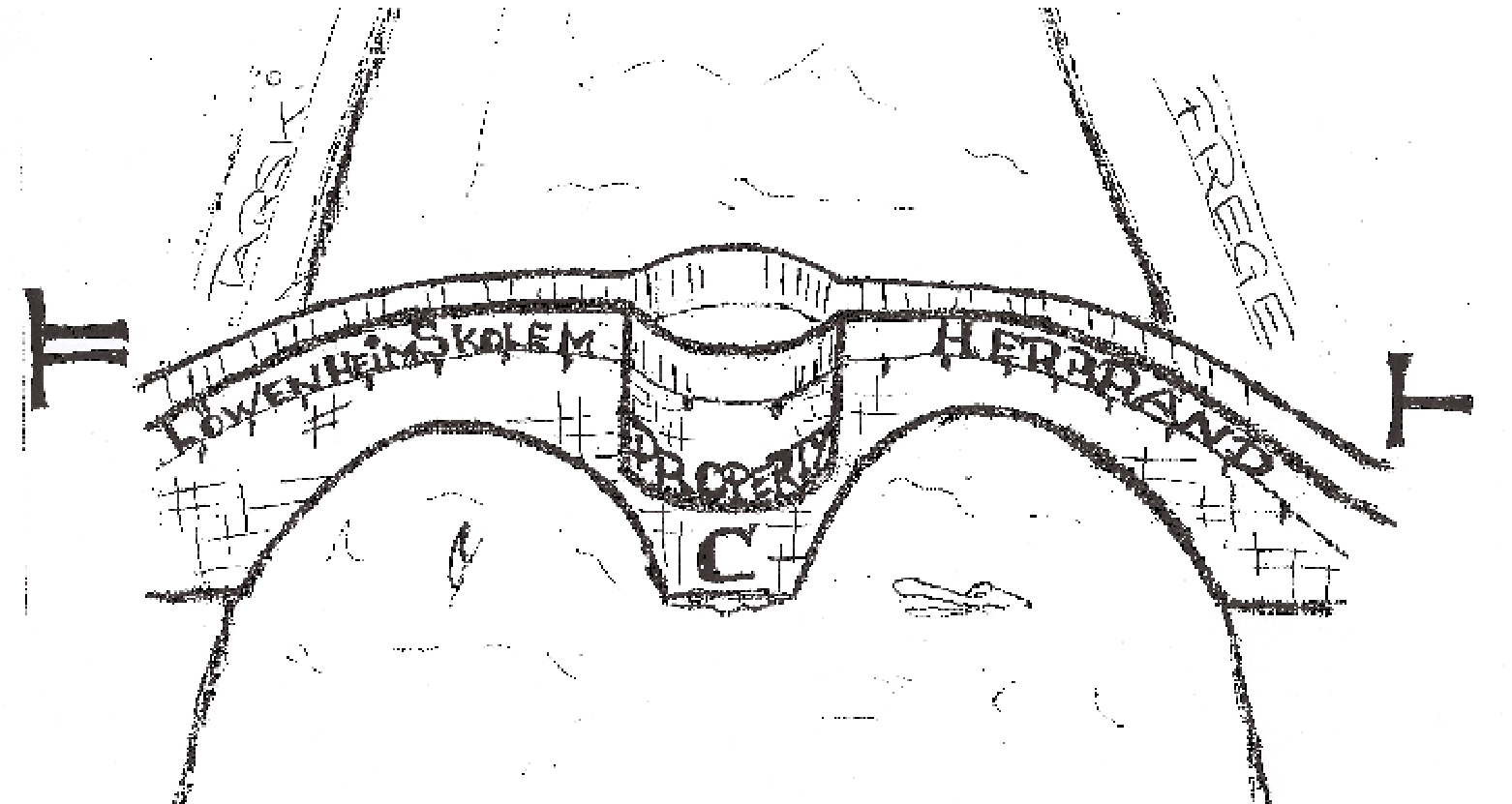}
\begin{minipage}{\textwidthminustwocm}{\footnotesize\mbox{}\\[+1.7ex]%
Figure: 
The bridge of the \protect\loewenheimskolemtheorem\ and 
\protect\herbrandsfundamentaltheorem, 
based on the sentential \protect\propertyC\ standing firm in the river 
that divides the banks of valid and derivable formulas in the land of 
first-order predicate logic.\par}\end{minipage}
\notop\end{center}\end{figure}

\subsection{Constructiveness of \herbrandsfundamentaltheorem}

Why was the difference between the model-theoretic notion of validity
and the constructive notion of derivability in a sound and complete calculus
so crucial for \herbrand?
The reason,
of course, 
is the undecidability of first-order logic,
which essentially requires the non-constructive
use of actual infinities in the definition of validity. \hskip.2em
\hilbert's program in \nolinebreak logic 
---~best described in \cite{grundlagen-german-english-edition-volume-one-one}~---
was to show the consistency of such non-constructive methods in mathematics
by {\em finitistic methods}, 
\ie\ by methods that are even more restrictive than the intuitionistic
methods in mathematics following \brouwername.

\herbrand\ does not accept any
  model-theoretic semantics unless the models are finite. \ In this
  respect, \herbrand\ is more finitistic than \hilbert, who demanded
  finitism only for consistency proofs. 
\notop\halftop\begin{quote}
``\herbrand's negative view of set theory leads him to take,
on certain questions,
a \nolinebreak stricter attitude than \hilbert\ and his collaborators. \ 
He is more royalist than the king. \
\hilbert's metamathematics has as its main goal to establish the consistency 
of certain branches of mathematics and thus to justify them; \hskip.3em
there, one had to restrict himself to finitistic methods. \ 
But in logical investigations other than the consistency problem of
mathematical theories the 
\index{Hilbert!school}%
\hilbert\ school was ready to work with
set-theoretic notions.''\getittotheright
{%
\index{Heijenoort!Jean van}%
\cite[\p 118]{heijenoort-work-herbrand}}\notop\end{quote}
As a consequence of this ``royalist'' attitude,
\herbrand\ was very proud on the fact that his \fundamentaltheorem\ 
is perfectly constructive in the sense that its proof shows how
anything claimed can be constructed from anything given: \hskip.3em
From \math A, we can construct an arbitrary large part of 
the sequence 
\maths{F^{\,\app{\termsofdepthnovars 1}F},
F^{\,\app{\termsofdepthnovars 2}F},
F^{\,\app{\termsofdepthnovars 3}F},
\ldots}. \hskip.6em
From a derivation of \maths A, 
we \nolinebreak can compute a number \math i such that 
\math A has \propertyC\ of order \nlbmath i
(\ie\ such that \math{F^{\,\app{\termsofdepthnovars i}F}} 
 is sententially valid).
\hskip.4em
If \math A has \propertyC\ of order \nlbmaths i, \hskip.2em
we can construct a linear derivation of \nlbmath A from 
\nlbmath{F^{\,\app{\termsofdepthnovars i}F}}
\ ---~~provided that we are explicitly 
given \nlbmath i as a definite number.%
\pagebreak

\section{Formal Presentation}
\subsection{Basic Notions and Notation}

Before we can present \herbrandsfundamentaltheorem\ formally,
we have to provide 
some further notions and notation on first-order formulas and
several inference rules for first-order logic.
Note that we will partly use modern notions,
which did not exist at \herbrand's time.

If we want to focus on a certain position in a formula,
we write the formula as \maths{A[B]}. \
This means that \math B is a formula that occurs in the 
context \math{A[\ldots]} as a sub-formula at a certain fixed position,
which, however, is not explicitly given by the notation.
Then we \nolinebreak denote with \nlbmath{A[C]} the formula that results
from the formula \math{A[B]} by replacing the one occurrence
of \nlbmath B at the fixed position with the 
formula \nlbmath C.

We denote with \bigmaths{A
\{x_1\tight\mapsto t_1,\ldots,x_n\tight\mapsto t_n\}}{}
the result of replacing all occurrences of the distinct variables
\nlbmath{x_1,\ldots,x_n} in the formula \nlbmath A in parallel with the
terms \nlbmaths{t_1,\ldots,t_n}, respectively. \hskip.3em
Here, \maths{\{x_1\tight\mapsto t_1,\ldots,x_n\tight\mapsto t_n\}}{}
\nolinebreak is a notation for a {\em substitution},
\ie\ for a function from variables to terms.

The occurrence of a quantifier in a formula is {\em accessible}\/
if it is not in the scope of any other quantifier.
For instance, 
in the valid formula 
\\\LINEmaths{\forall x\stopq\exists y\stopq\inpit{x\tightprec y}
\oder\exists m\stopq\forall z\stopq\neg\inpit{m\tightprec z}}{}
\\on the binary predicate symbol \nlbmath\prec\ (with infix notation), \hskip.2em
the occurrences of the quantifiers 
\math{\forall x.} \nolinebreak and \nlbmath{\exists m.}
are the only accessible ones. \hskip.3em
Note that we assume the scopes of our quantifiers to be minimal
in the sense that the scope of \hskip.1em\math{\forall x.} in this formula
does not 
include the sub-formula \bigmathnlb{
\exists m\stopq\forall z\stopq\neg\inpit{m\tightprec z}}{}
---~~contrary to the formula 
\\\LINEmaths{\forall x\stopq\inparenthesestight{\exists y
\stopq\inpit{x\tightprec y}\oder
\exists m\stopq\forall z\stopq\neg\inpit{m\tightprec z}}},
\\where only the occurrence of \nolinebreak\hskip.1em\math{\forall x.} 
is accessible.

\citet{smullyan} classified reductive inference rules
---~and the inference rules of the \hilbert\ calculi 
we will consider here can all be seen as such if we read them
\mbox{bottom up~---}
into
\math\alpha~(sentential+non-branching),
\math\beta~(sentential+\underline branching),
\math\gamma, 
and \nlbmath\delta. \hskip.2em
According to this classification,
we introduce the following notion on quantifiers,
bearing in mind that \maths\tightund, 
\nlbmaths\tightoder, \nolinebreak and \nlbmath\neg\ \hskip.1em
are our only \myboolean\ operators.

The occurrence of a quantifier in formula is {\em\existentialoid} \hskip.1em
if it is of the form \nlbmath{\exists x.} 
and it is in the scope of an even number of
negation symbols,
or of the form \nlbmath{\forall x.} 
and in the scope
of an odd number of negation symbols;
otherwise the quantifier is {\em\universaloid}. \hskip.3em
(\Anexistentialoidquantifier\ turns up as \nlbmath\exists\ in a prenex form of
 the formula, 
 and \auniversaloidquantifier\ as \nlbmaths\forall.)

The occurrence of a variable in a formula is {\em\existentialoid}\/ \hskip.1em
if it is bound by \anexistentialoidquantifier; \hskip.2em
it is \universaloid\ if \nolinebreak it \nolinebreak is 
bound by \auniversaloidquantifier\ or free 
(\ie\ not bound by any quantifier).

\vfill\pagebreak
\subsection{A Modern Version of \herbrand's {\it Modus Ponens}\/-Free Calculus}

Now we are prepared to understand the following three inference rules
which constitute a slightly improved version of 
\herbrand's {\it modus ponens}\/-free calculus in the style of 
\makeacitetofthree
{heijenoort-herbrand}
{heijenoort-modern-logic}
{heijenoort-work-herbrand}
and \makeacitetoftwo
{wirth-heijenoort}
{wirth-heijenoort-SEKI}.

Note that we may rename bound variables to satisfy the side conditions
of the inference rules, 
because we consider equality of formulas only up to renaming of bound
variables.

\par\yestop\noindent{\bf\Ruleofexistentialoidquantification:} \ 
\bigmaths{\begin{array}[c]{l}{A[H\{x\mapsto t\}]}
\\\hline
  {A[Q x\stopq H]}
\\\end{array}}{~}
where\begin{enumerate}\item
\math{Q x.} is an accessible \existentialoidquantifier\ of 
\nlbmath{A[Q x\stopq H]}, \ and\item the free variables of the term \nlbmath t 
\hskip.1em must not be bound by quantifiers in \nlbmaths H.%
\yestop\end{enumerate}

\begin{example}[Application of the \ruleofexistentialoidquantification]\label
{example Application of the rule of existentialoid quantification}%
\mbox{}
\\\noindent If the variable \math z does not occur free in the term \nlbmaths t,
we get the following two inference steps with identical premises
by application of the \ruleofexistentialoidquantification\
at two different positions:\begin{itemize}\item
\math{\begin{array}[c]{r c r}\inpit{t\tightprec t}
 &\oder
 &\neg\forall z\stopq\inpit{t\tightprec z}
\\\hline
  \inpit{t\tightprec t}
 &\oder
 &\exists x\stopq\neg\forall z\stopq\inpit{x\tightprec z}
\\\end{array}}{} \
via the meta-level substitution
\getittotheright{\maths{
\{\ \ \ \ A[\ldots]\ 
\mapsto\ \inpit{t\tightprec t}
\oder [\ldots]\comma\ \ \
H\ \mapsto\ \neg\forall z\stopq\inpit{x\tightprec z}
\comma\ \ \
Q\ \mapsto\ \exists
\ \ \ \ \}\majorfootroom\majorheadroom};\mbox{~~~~~\,}}\item
\maths{\begin{array}[c]{r c r}\inpit{t\tightprec t}
 &\oder
 &\neg\forall z\stopq\inpit{t\tightprec z}
\\\hline
  \inpit{t\tightprec t}
 &\oder
 &\neg\forall x\stopq\forall z\stopq\inpit{x\tightprec z}
\\\end{array}}{} \
via the meta-level substitution
\getittotheright{\maths{
\{\ \ \ \ A[\ldots]\ \mapsto\ 
\inpit{t\tightprec t}
\oder\,\neg[\ldots]
\comma\ \ \
H\ \mapsto\ \forall z\stopq\inpit{x\tightprec z}
\comma\ \ \
Q\ \mapsto\ \forall
\ \ \ \ \}\majorfootroom\majorheadroom}.}\qed\end{itemize}
\end{example}

\par\yestop\noindent
{\bf\Ruleofuniversaloidquantification:} \bigmaths{\begin{array}[c]{l}{A[H]}
  \\\hline
    {A[Q y\stopq H]}
  \\\end{array}}{~}
where\begin{enumerate}\item 
\math{Q y.} is an accessible \universaloidquantifier\
of \nlbmaths{A[Q y\stopq H]}, \ and\item
the variable \nlbmath y must not occur free 
in the context \nlbmath{A[\ldots]}.\majorfootroom\end{enumerate}

\par\yestop\noindent
{\bf\Ruleofsimplification:} \bigmaths{\begin{array}[c]{l}{A[H\circ H']}
\\\hline
  {A[H]}
\\\end{array}}{~} where\begin{enumerate}\item
``\math{\circ}'' stands for ``\tightoder'' if \math{[\ldots]} \nolinebreak
occurs in the scope of an even number of negation symbols 
in \nlbmath{A[\ldots]}, \ 
and for ``\tightund'' otherwise, 
and\item
\math{H'} is a variant of the sub-formula \nlbmath H \
(\ie, \math{H'} is \math H or 
 can be obtained from \nlbmath H by the renaming of variables
bound in \nlbmath H).\majorfootroom\end{enumerate}
Moreover, the 
{\em generalized rule of \math\gamma-simplification} 
  is the sub-rule for the case that
  \math H \nolinebreak is of the form \bigmaths{Q y\stopq C}{}
  and \bigmaths{Q y.}{} is \anexistentialoidquantifier\
  of \bigmaths{A[Q y\stopq C]}.%
\pagebreak\par
\begin{remark}%
[Historic Version of \herbrand's {\it Modus Ponens}-Free Calculus]%
\label{remark historically correct calculus}%
\par\noindent\sloppy
The before-mentioned three rules 
are to be used for a modern presentation of \herbrand's
{\it modus ponens}\/-free calculus.
The historical {\it modus ponens}\/-free calculus of \herbrand\ 
had the generalized rule of simplification,
but
only the shallow rules of 
\mbox{``\math\gamma- and \math\delta-quantification\closequotecommanospace} \
compensated by the addition of the rules of passage.

\halftop\noindent
{\bf Rules of \math\gamma- and \math\delta-quantification}
result from our formalization of the generalized rules 
by restricting \nlbmath{A[\ldots]} to the empty context
\hskip.1em
(\ie\nolinebreak\ \math{A[Q x\stopq H]}, \eg, is just 
\nolinebreak\mbox{~\math{Q x\stopq H}).} 

\par\halftop\noindent{\bf Rules of Passage:} \
The following six \nolinebreak logical
  equivalences may be used for rewriting from left to right 
({\em prenex direction}\/) and from right to left 
({\em anti-prenex direction}\/), 
resulting in twelve \nolinebreak deep 
inference rules
(where \math B is a formula 
 in which the variable \nlbmath x does not occur free):
\\\noindent\LINEmath{\begin{array}{l c r@{~~~~}c@{~~~~}l}
       (1)
      & 
      &\neg\forall x\stopq A
      &\equivalent
      &\exists x\stopq \neg A
     \\(2)
      & 
      &\neg\exists x\stopq A
      &\equivalent
      &\forall x\stopq \neg A
     \\(3)
      & 
      &\inpit{\forall x\stopq A}\nottight\oder B
      &\equivalent
      &\forall x\stopq \inpit{A\hskip.06em\tightoder B}
     \\(4)
      & 
      &B\nottight\oder\forall x\stopq A
      &\equivalent
      &\forall x\stopq \inpit{B\hskip.06em\tightoder A}
     \\(5)
      & 
      &\inpit{\exists x\stopq A}\nottight\oder B
      &\equivalent
      &\exists x\stopq \inpit{A\hskip.06em\tightoder B}
     \\(6)
      & 
      &B\nottight\oder\exists x\stopq A
      &\equivalent
      &\exists x\stopq \inpit{B\hskip.06em\tightoder A}
      \\\end{array}}%
\par\halftop\noindent
Note that \herbrand\ did not need rules of passage for conjunction
(besides the rules of passage for negation (1,\,2) and for disjunction 
(3,\,4,\,5,\,6)),
\hskip.2em
because he considered conjunction 
\bigmaths{A\hskip.09em\tightund\hskip.05em B}{} a
meta-level notion defined as \bigmaths
{\neg\inpit{\neg A\hskip.14em\tightoder\hskip.10em\neg B}}.

\herbrand\ needed his rules of passage
(in anti-prenex direction)
for the completeness of his historic {\it modus ponens}\/-free calculus
because the shallow rules of quantification
\mbox{---~contrary} to the generalized \mbox{ones~---}
cannot introduce quantifiers at 
non-top positions.

\herbrandindex\herbrand\ introduced these rules 
in \litsectref{2.2} of his \PhDthesis\
\cite{herbrand-PhD}. \
He \nolinebreak named 
the rules of \mbox{\math\gamma- and \math\delta-quantification}
``second'' and ``first rule of generalization'' 
\cite[\p\,74\f]{herbrand-logical-writings},
respectively 
({\it ``deuxi\`eme''}
 and 
 {\it ``premi\`ere r\`egle de g\'en\'eralisation''}
 \cite[\p\,68\f]{herbrand-ecrits-logiques}). \hskip.4em
At the same places, 
we also find 
the ``rules of passage'' 
({\it\mbox{``r\`egles} de passage''}\/\nolinebreak\hskip.2em\nolinebreak).
\hskip.4em
Finally, in \litsectref{5.6.A} of his \PhDThesis,
\herbrandindex\herbrand\ also introduces the generalized rule of simplification
\cite[\p\,175]{herbrand-logical-writings} 
({\it``r\`egle de simplification g\'en\'eralis\'ee''}
\cite[\p\,143]{herbrand-ecrits-logiques}).
\getittotheright\qed\end{remark}

\vfill\pagebreak
\subsection{\propertyC}
\halftop\begin{definition}[Height of a Term, Champ Fini 
\nlbmath{\termsofdepth n F}]\label
{definition champs finis}\mbox{}
\par\noindent
We \nolinebreak use \nlbmath{\CARD t}
to denote the {\em height}\/ of a term \nlbmath t,
which is given by
\par\noindent\LINEmaths{\CARD{\anonymousfpp{t_1}{t_{m}}}
  \nottight{\nottight{\nottight=}}1+\max\{0,\CARD{t_1},\ldots,\CARD{t_m}\}}.
\par\noindent For a positive natural number \nlbmath n
and a formula \nlbmath F,
as a finite substitute for a typically infinite, 
full term universe,
\herbrandindex\herbrand\ uses what he calls a 
{\em champ fini of order \nlbmaths n,} \hskip.1em
which we \nolinebreak will denote with 
\nlbmath{\termsofdepth n F}.  \ 
The terms of \nlbmath{\termsofdepth n F} 
are constructed from the symbols that occur free
in \nlbmath F: \
the function symbols, the constant symbols
(which we will tacitly subsume 
 under the function symbols in what follows), \hskip.1em
and the free variable symbols
(which can be seen as constant symbols here). \hskip.3em
Such a 
{\em champ fini}\/ differs from a full term universe in
containing only the terms \nlbmath t with \bigmaths{\CARD t\prec n}{\,.} \ 
\par\noindent 
So we have \bigmaths{\termsofdepth 1 F\tightequal\emptyset}. \hskip.2em
\par\noindent To guarantee \bigmaths
{\termsofdepth n F\tightnotequal\emptyset}{} \hskip.1em
for \bigmaths{n\succ 1},
in case that neither constants nor free variable symbols occur in \nlbmaths F,
\hskip.1em we will assume that a fresh constant symbol 
\nolinebreak ``\nlbmath\bullet\nolinebreak\hskip.02em\nolinebreak''
(which does not occur elsewhere)
\hskip.1em is \nolinebreak
included in the term construction in addition to the free symbols of \nlbmath F.
\getittotheright\qed\end{definition}

\halftop\halftop\noindent\herbrand's definition of an expansion 
follows the traditional
idea that ---~for a {\em finite}\/ domain~--- universal (existential) 
quantification can be seen as a finite conjunction (disjunction) over the 
elements of the domain:
\begin{definition}[Expansion]\label
{definition sub-expansion}\sloppy\mbox{}\\Let \nlbmath{\mathcal T} be a 
finite set of terms. \
To simplify substitution,
let \nlbmath A be a formula whose bound variables do not occur in 
\nlbmaths{\mathcal T}.
\\The 
{\em expansion \nlbmath{A^{\mathcal T}} 
of \nlbmath A \wrt\ \nlbmath{\mathcal T}} \hskip.3em
is the formula given by the following recursive definition.
\\If \math A is quantifier-free formula,
then \bigmaths{A^{\mathcal T}:=A}.
Moreover:
\bigmaths{\inpit{\neg A_1}^{\mathcal T}:=\neg A_1^{\mathcal T}},
\\\bigmaths{\begin{array}[b]{l l l}\inpit{A_1\oder A_2}^{\mathcal T}
 &:=
 &A_1^{\mathcal T}\oder A_2^{\mathcal T},
\\\inpit{A_1\und A_2}^{\mathcal T}
 &:=
 &A_1^{\mathcal T}\und A_2^{\mathcal T},
\\\end{array}}{}
\bigmaths{\begin{array}[b]{l l l}\inpit{\exists x.\,A}^{\mathcal T}
 &:=
 &\bigvee_{t\in\mathcal T}\,A^{\mathcal T}\{x\tight\mapsto t\},
\\\inpit{\forall x.\,A}^{\mathcal T}
 &:=
 &\bigwedge_{t\in\mathcal T}\,A^{\mathcal T}\{x\tight\mapsto t\}.
\\\end{array}}{}
\getittotheright\qed\end{definition}
\begin{definition}[Outer \skolemizedForm]\label
{definition skolemized forms}\mbox{}\\
The {\em outer \skolemizedform\ of a formula \nlbmath A}\/ results from \nlbmath
A by removing every \universaloidquantifier\ and replacing its bound
variable \nlbmath x with \nlbmaths{\app{\forallvari x{}}{y_1,\ldots,y_m}}, \ 
where \math{\forallvari x{}} \nolinebreak is a fresh (``\skolem'') symbol
and \maths{y_1,\ldots,y_m}, in this order,
are the variables of the \math\gamma-quantifiers in
whose scope the \mbox{\math\delta-quantifier} occurs.
\getittotheright\qed\end{definition}

  
\begin{definition}[\propertyC]\label
{definition properties C and C star}\mbox{}
\\Let \math A be a \firstorder\ formula. \ 
Let \math n be a positive natural number. \
\\Let \math F be
the outer \skolemizedform\ of \nlbmaths A.
\\\math A \nolinebreak
{\em has \propertyC\ of order \nlbmath 1}\/ \hskip.2em if \hskip.2em
\math F is a \sententialtautology.
\\For \bigmaths{n>1}, the formula \math A \nolinebreak
{\em has \propertyC\ of order \nlbmath n}\/ \hskip.2em if\\the 
expansion \nlbmath{F^{\,\app{\termsofdepthnovars n}F}} is a
\sententialtautology.
\getittotheright\qed\end{definition}

\vfill\pagebreak

\subsection{The Theorem and its Lemmas}
\halftop\halftop
\begin{theorem}[\herbrandsfundamentaltheorem\ \`a la \heijenoort]%
\label{theorem herbrand fundamental two}%
\\Let\/ \math A be a \firstorder\ formula. 
The following two statements are logically equivalent. \hskip.3em
Moreover, 
we can construct a witness for each statement 
from a witness for the other one.\begin{enumerate}\notop\item[1.]%
There is a positive natural number \nlbmath n such that
\math A has \propertyC\ of order \nlbmaths n.\noitem\item[2.]\sloppy 
There\,is\,a\,\sententialtautology\ \nlbmaths B, 
and 
\\there\,is\,a\,derivation\,of \math A from \nlbmath B
that\,consists\,in\,applications\,of the generalized rules of
\\simplification, \hskip.2em
\mbox{\math\delta-quantification,} \hskip.2em
and\/ \mbox{\math\gamma-quantification} \hskip.2em
\\(and in the renaming of bound variables).
\getittotheright\qed\end{enumerate}\end{theorem}
\par\halftop\halftop\noindent
As we can decide \propertyC\ of order \nlbmath n \hskip.1em for 
\maths{n\tightequal 1}, \maths{n\tightequal 2},
\maths{n\tightequal 3}, \maths\ldots, \hskip.1em
\theoref{theorem herbrand fundamental two}
immediately provides us with a semi-decision procedure for derivability
(and, thus, by the \loewenheimskolemtheorem, also for validity)
of any \firstorder\ formula \nlbmath A given as input.

Note that the witnesses mentioned 
in \theoref{theorem herbrand fundamental two}
are, of course,
on the one hand, 
a \nolinebreak concrete representation of 
the natural number \nlbmaths n, \hskip.1em
and, on the other hand, concrete representations of 
the formula \nlbmath B and of
the derivation of \math A from \nlbmaths B.

To get some more information on the construction of these witnesses,
we have to decompose the equivalence of \theoref
{theorem herbrand fundamental two}
into the two implications found in the following two lemmas,
which constitute the theorem.

\begin{lemma}
[From \propertyC\ to a Linear Derivation]
\label{lemma from C to yields a la heijenoort}%
\\\noindent
Let\/ \math A be a \firstorder\ formula. \
Let\/ \math F be the outer\/ \skolemizedform\ of \nlbmath A. \
Let\/ \math n be a positive natural number.
\\\noindent If\/ \math A has \propertyC\ of order \nlbmath n, 
then we can construct a derivation of \nlbmath A
of the following form, 
in which we read any term starting with a\/ \skolem\ function 
as an atomic variable:%
\\\noindent\begin{tabular}[b]{@{}l l}%
  \headroom
  {\bf Step\,1:} 
 &We \nolinebreak start with 
  the \sententialtautology\/ \nlbmaths{F^{\,\termsofdepth n F}\!}.
\\\headroom{\bf Step\,2:}
 &Then we may repeatedly apply the generalized rules of\/
  \math\delta- and \math\gamma-quanti\-fi\-cation.
\\\headroom{\bf Step\,3:}
 &Then we may repeatedly apply the generalized rule of\/
  \math\gamma-simplification.
\\\headroom{\bf Step\,4:}
 &Then we rename all bound \mbox{\math\delta-variables} 
  to obtain \nlbmath A.\\\end{tabular}\\[-2.5ex]\getittotheright\qed\end{lemma}
\par\halftop\halftop\noindent
The proof idea of \lemmref{lemma from C to yields a la heijenoort}
is to transform the computation of the expansion of the 
outer \skolemizedform\ into a 
reduction in \herbrand's {\it modus ponens}\/-free calculus.
In this transformation,
the reduction with the generalized rules of \math\gamma-simplification
and \mbox{\math\gamma-quantification} models the expansion,
and the renaming of bound \math\delta-variables 
to \skolem\ terms considered as variable names
models the \skolemization,
whereas the reduction with the generalized rules of
\math\delta-quantification just drops the \math\delta-quantifiers.
The critical task in this transformation is to find an appropriate total
order of the variable occurrences of the original expansion steps,
so that the side conditions of 
the resulting reductive applications of the inference rules are met. 

See \nolinebreak\cite[\litsectref 5]{wirth-heijenoort-SEKI}
for an elaborate, but easily conceivable example
for an application of a procedure that can actual construct
such a derivation.
That example also shows how to overcome the inefficiency of this procedure
and how to find a proof of a manageable size.

\pagebreak

\begin{lemma}[From a Linear Derivation to \propertyC]
\label{lemma from yields to C a la heijenoort}\sloppy\\\noindent
If there is a derivation of the \firstorder\ formula \nlbmath A 
from a \sententialtautology\
by applications of the 
generalized rules of simplification, and of\/
\math\gamma- and\/ \math\delta-quantification 
\mbox{(and renaming of bound variables),}
\\then \math A has \propertyC\ of order\/ 
\nlbmaths{\displaystyle1+\sum_{i=1}^m\CARD{t_i}},
\par\noindent
where \math{t_1,\ldots,t_m} are the instances for the meta-variable \nlbmath t
of the generalized rule of\/ \mbox{\math\gamma-quantification} in its \math m
applications in the derivation of \nlbmaths A.
\getittotheright\qed\end{lemma}

\halftop
\begin{remark}[Historical Version of \herbrandsfundamentaltheorem]
\par\noindent
As already explained in \remaref{remark historically correct calculus},
\herbrand's actual calculus was a bit different and had to take 
the detour via adding quantifiers on top level and then moving them in.
This seemed to admit a minor simplification by a 
detour via the prenex normal form.
To reduce a problem to problems of manageable size
({\it divide et impera}\/), \hskip.2em
the detour via prenex normal form was a leading standard 
at \herbrand's time. 
Meanwhile prenex normal form plays a lesser \role\
in the better logic courses
because of its crucial efficiency problems.
\par
In \nolinebreak \herbrand's case this problem turned out to be fatal for the
correctness of his proof:
\herbrand\ computed the upper bound for the order of 
\propertyC\ after application of the rules of passage much lower than
it actually is. This 
is well-documented under the name of \herbrandsfalselemma. \hskip.5em
\par
One correction of \herbrandsfalselemma\ 
is the one that we have presented in this article
and that consists in
adding 
---~to \herbrand's deep version of his inference rule of simplification~---
also the deep versions of his inference rules of quantification. \hskip.3em
Looking at the style in which the great mathematician \herbrandname\
organized his most creative work in logic we may say that,
if anybody had noticed this bug 
in \herbrand's proof during \herbrand's lifetime,
this correction would have been the most straightforward \mbox{bug fix} for him.
\hskip.3em
Moreover, 
this correction still is 
the most straightforward and most elegant one today. \hskip.3em
It was clearly outlined by \heijenoortname, \hskip.1em
but first sketched in publication in \cite{herbrand-handbook}, \hskip.3em
and first published with an explicit presentation 
in \cite{wirth-heijenoort}.\end{remark}

\halftop\halftop\halftop\halftop\section{Conclusion}
In this article
we have delivered what we consider the very essentials that any logician
should know on \herbrandsfundamentaltheorem, and we suggest 
\cite{wirth-heijenoort-SEKI} and \cite{SR--2009--01} for further reading on 
\herbrandsfundamentaltheorem, \herbrand's further work in logic,
and for a listing of further sources on the subject.
\vfill\pagebreak

\small

\addcontentsline{toc}{section}{Acknowledgments}
\section*{Acknowledgments}
I would like to thank \anellisname\ \anellislifetime\ for all his kind help
for me and my work on \herbrandname\ and \heijenoortname.
\irving\ has made my life better and richer.
His death has left a gap in my life.

\nocite{writing-mathematics}
\addcontentsline{toc}{section}{References}
\bibliography{herbrandbib}

\end{document}